\documentstyle[12pt]{article}
\newcommand{\sect}[1]{\section{#1}\setcounter{equation}{0}}

\font\mbn=msbm10 scaled \magstep1 \font\mbs=msbm7 scaled \magstep1
\font\mbss=msbm5 scaled \magstep1
\newfam\mbff
\textfont\mbff=\mbn \scriptfont\mbff=\mbs
\scriptscriptfont\mbff=\mbss\def\mbf{\fam\mbff}
\def\Re{{\mbf R}}

\def\Z{{\mbf Z}}
\def\Co{{\mbf C}}

\def\Di{{\mbf D}}
\def\Bo{{\mbf B}}
\def\N{{\mbf N}}

\newtheorem{Th}{Theorem}[section]
\newtheorem{Lm}[Th]{Lemma}
\newtheorem{C}[Th]{Corollary}

\newtheorem{D}[Th]{Definition}
\newtheorem{Proposition}[Th]{Proposition}
\newtheorem{R}[Th]{Remark}
\newtheorem{E}[Th]{Example}
\author{Alexander Brudnyi\thanks{Research supported in part by NSERC.
}\\
Department of Mathematics and Statistics\\
University of Calgary, Calgary\\
Canada}
\title{On Local Behavior of Holomorphic Functions Along Complex Submanifolds of $\textbf{C}^{N}$}
\date{}
\begin{document}
\maketitle
\begin{abstract}
{In this paper we establish some general results on local behavior of holomorphic functions along complex submanifolds of $\Co^{N}$. As a corollary, we present multi-dimensional generalizations of an important result of
Coman and Poletsky on Bernstein type inequalities on transcendental curves in $\Co^{2}$.}
\end{abstract}
\sect{\hspace*{-1em}. Formulation of the Main Result} 
{\bf 1.1.} In this paper we establish some general results on restrictions of holomorphic functions to complex submanifolds of $\Co^{N}$.  The subject pertains to the area of the, so-called, polynomial inequalities for analytic and plurisubharmonic functions that includes, in particular, Bernstein, Markov and Remez type inequalities.
Recently there has been a
considerable interest in such inequalities in connection with various problems of analysis. Let us recall that the classical univariate
inequalities for polynomials have appeared in approximation theory and for a long time have been considered as technical tools for proofs of Bernstein type inverse theorems. At the present time polynomial type inequalities have been found a lot of important applications in areas which are well apart from approximation theory. We will only briefly mention several of these areas.

The papers [GM], [Bou] and [KLS] apply polynomial inequalities with different integral norms to study some problems of Convex Geometry (in particular, the famous Slice Problem).

In the papers [B1], [B2], [BB], [G], [P] and [PP] and books [DS] and [JW]
Chebyshev-Bernstein and related Markov type inequalities are used to explore a wide range of properties of the classical spaces of smooth
functions including Sobolev type embeddings and trace theorems, extensions and differentiability.

The papers [FN1] and [FN2] on Bernstein type inequalities for traces of polynomials to algebraic varieties were inspired by and would have important applications to some basic problems of the theory of subelliptic differential equations.

The paper [BLMT] discovers a profound relation between the 
exponents in the tangential Markov inequalities for restrictions of polynomials to a smooth manifold $M\subset\Re^{N}$ and the property of $M$ to be an algebraic manifold.

An application of polynomial inequalities to Cartwright type theorems for entire functions is presented in [Br1]
and [Br2], see also [LL], [Lo], [K].

In [T1], [T2] Bernstein type inequalities are used to obtain
new results in transcendental number theory.

Finally, we mention applications of polynomial inequalities to the second part of Hilbert's sixteenth problem concerning the number of limit cycles of planar polynomial vector fields, see [I], [RY], [Br3] and [Br4].

In [CP] Coman and Poletsky obtained an important result on Bernstein type inequalities for restrictions of holomorphic polynomials to certain transcendental curves in $\Co^{2}$. The main purpose of our paper is to present a general approach to such kind of inequalities. As an application, we obtain  multi-dimensional generalizations of the result of [CP].
\\
\\
{\bf 1.2.} In what follows by $\Bo_{r}^{n}(z_{0})\subset\Co^{n}$ we denote the open Euclidean ball of radius $r$ centered at $z_{0}$. We set 
$$
\Bo_{r}^{n}:=\Bo_{r}^{n}(0),\ \ \ \Bo^{n}:=\Bo_{1}^{n}, \ \ \
\Di_{r}(z_{0}):=\Bo_{r}^{1}(z_{0}),\ \ \
\Di_{r}:=\Bo_{r}^{1},\ \ \ \Di:=\Bo^{1}.
$$
By $\overline{S}$ and $\partial S$ we denote the closure and the boundary of $S\subset\Co^{n}$.

For a continuous function
$f:\Bo_{r}^{n}(z_{0})\to\Co$ we define
$$
M_{f}(r,z_{0}):=\sup_{\Bo_{r}^{n}(z_{0})}|f|,\ \ \  m_{f}(r,z_{0}):=\ln M_{f}(r,z_{0}).
$$
If $z_{0}=0$ we set $M_{f}(r):=M_{f}(r,0)$, $m_{f}(r):=m_{f}(r,0)$.

Next, by ${\cal L}_{n}$ we denote the family of complex lines $l\subset\Co^{n}$ passing through the origin. For each $l\in {\cal L}_{n}$ we naturally identify
$l\cap\Bo_{r}^{n}$ with $\Di_{r}$.

The main result of the paper (Theorem \ref{te1}) consists of two inequalities that express in a quantitative way the following simple fact.

{\em Let $f$ be a nonconstant holomorphic function in $\Bo_{r}^{n}$ and $g$ be a holomorphic function in the domain $\Omega:=\Bo_{r}^{n}\times\Di_{M_{f}(r)}\subset\Co^{n+1}$. 
Suppose that 
\begin{itemize}
\item[(1)]
There exist a complex line $l\in {\cal L}_{n}$, a vector
$(0,v)\in\Co^{n+1}$, $v\in\Co$, and a positive number $N\in\N$ such that the affine line
$l(v):=l+(0,v)$ intersects the graph $\Gamma_{f}:=\{(z,f(z))\in\Co^{n+1}\ :\ z\in\Bo_{r}^{n}\}$ transversely in at least $N$ points;
\item[(2)]
The univariate holomorphic function $g|_{l(v)\cap\Omega}$ is nonconstant and has less than $N$ zeros.
\end{itemize}
Then the restriction $g|_{\Gamma_{f}}$ is not identically zero.}

As a corollary of Theorem \ref{te1} we obtain Bernstein, Markov, Remez and Jensen type inequalities for $g|_{\Gamma_{f}}$, see section 2.1. Also,
for $f$ being a  nonpolynomial entire function on $\Co^{n}$ our result leads to a generalization of one of the main results of [CP], see section 2.2.

To give the precise formulation of the result we introduce some notations.

Assume that $f:\Di_{r}(z_{0})\to\Co$ is holomorphic. By $n_{f}(r,z_{0})$ we denote the number of zeros of $f$ in $\Di_{r}(z_{0})$. (We write $n_{f}(r,z_{0})=-\infty$ if $f\equiv 0$.) Then the valency of $f$ in $\Di_{r}(z_{0})$ is defined by 
$$
v_{f}(r,z_{0}):=\sup_{c\in\Co}n_{f+c}(r,z_{0}).
$$ 

Also, the Bernstein index $b_{f}$ of $f$ is given by the formula
$$
b_{f}(r,z_{0}):=\sup\{m_{f}(es,z)-m_{f}(s,z)\}
$$
where the supremum is taken over all
$\Di_{es}(z)\subset\subset\Di_{r}(z_{0})$. (If $f\equiv 0$ we assume that $b_{f}(r,z_{0})=0$.) 

Let us mention that the values of $n_{f}$, $v_{f}$ and $b_{f}$ are finite for a nonzero $f$ defined in a neighbourhood of the closure of $\Di_{r}(z_{0})$. We set for brevity
$$
n_{f}(r):=n_{f}(r,0),\ \ \ v_{f}(r):=v_{f}(r,0),\ \ \
b_{f}(r):=b_{f}(r,0).
$$
{\bf 1.3.} We are ready to formulate the main result of the
paper.

Suppose that $f$ is a holomorphic function in $\Bo_{tr}^{n}$, $r>0$, $1<t\leq 9$, satisfying 
\begin{equation}\label{e12}
M_{f}(r/t)\geq M_{1},\ \ \ M_{f}(tr)\leq M_{2}\ \ \ {\rm and}\ \ \
R_{f}(r,t,t^{2})\geq t
\end{equation}
where 
\begin{equation}\label{R}
R_{f}(r,t,s):=\frac{M_{f}(r/t)}{M_{f}(r/s)}\ \!,\ \ \
t\leq s<\infty.
\end{equation}
(The last inequality in (\ref{e12}) is valid, e.g., if $f(0)=0$.)

For every $l\in {\cal L}_{n}$ we set
$$
f_{l}:=f|_{l\cap\Bo_{tr}^{n}}
$$
and determine positive numbers $V_{f}(r,t)$ and $N_{f}(r,t)$ by the formulas
\begin{equation}\label{e10}
V_{f}(r,t):=\inf_{l}\{v_{f_{l}}(r/\sqrt{t})\ :\ f_{l}\neq const\},
\end{equation}
\begin{equation}\label{e11}
N_{f}(r,t):=\left\{
\begin{array}{ccc}
v_{f}(r/\sqrt{t})&{\rm if}&n=1\\
\\
\displaystyle\max\left\{\sup_{s\in [t,\infty)}\left\{\frac{\ln(R_{f}(r,t,s)/\sqrt{t})}{k(t,s)}\right\}, V_{f}(r,t)\right\}&{\rm if}&n\geq 2
\end{array}
\right.
\end{equation}
where
\begin{equation}\label{k}
k(t,s):=\ln\left(\frac{8e^{\pi^{2}}s\sqrt{t}}{(\sqrt{t}-1)^{2}}\right).
\end{equation}

In Lemma \ref{le41} we show that if $n=1$, then
$$
\frac{\ln R_{f}(r,t,s)}{k(t,s)}\leq v_{f}(r/\sqrt{t})\ \ \
{\rm for\ each}\ \ \ s\in [t,\infty).
$$
This, in particular, implies that
$$
\inf_{l}v_{f_{l}}(r/\sqrt{t})\leq N_{f}(r,t)\leq \sup_{l}v_{f_{l}}(r/\sqrt{t}).
$$

Let $g$ be a holomorphic function in the domain
$\Bo_{tr}^{n}\times\Di_{3M_{2}}\subset\Co^{n+1}$. For every $l\in {\cal L}_{n}$ we determine
\begin{equation}\label{e01}
g_{l}:=g|_{\Omega_{l}}\ \ \ {\rm where}\ \ \ \Omega_{l}:=(l\cap\Bo_{tr}^{n})\times\Di_{3M_{2}}.
\end{equation}
\begin{D}\label{d11}
We say that the function $g$ belongs to the class
${\cal F}_{p,q}(r;t;M_{2})$ for some $p,q\geq 0$ if 
\begin{equation}\label{e13}
\begin{array}{c}
M_{g_{l}(\cdot,w)}(tr)\leq e^{p}\cdot M_{g_{l}(\cdot,w)}(r)\ \ \
{\rm for\ all}\ \ \ l\in {\cal L}_{n},\ w\in\Di_{3M_{2}}
\ \ \ {\rm and}\\
\\
b_{g(z,\cdot)}(3M_{2})\leq q\ \ \ {\rm for\ all}\ \ \
z\in\Bo_{tr}^{n}.
\end{array}
\end{equation}
\end{D}

Set
\begin{equation}\label{e14}
g_{f}(z):=g(z,f(z)),\ \ \ \ z\in\Bo_{tr}^{n}.
\end{equation}
(This function is the restriction of $g$ to the graph $\Gamma_{f}\subset\Co^{n+1}$ of $f$.)

Our main result is
\begin{Th}\label{te1}
Assume that 
\begin{equation}\label{cond1}
p\leq \ln\left(\frac{1+t}{2\sqrt{t}}\right)\cdot N_{f}(r,t).
\end{equation}
Then there are positive constants $a_{1}(t), a_{2}(t)$ such that
for any $g\in {\cal F}_{p,q}(r;t;M_{2})$
\begin{equation}\label{e15}
\sup_{\Bo_{r}^{n}\times\Di_{M_{2}}}|g|\leq
\left(\frac{a_{1}(t)M_{2}}{M_{1}}\right)^{a_{2}(t)(p+q)}M_{g_{f}}(r),
\end{equation}
\begin{equation}\label{e18}
M_{g_{f}}(tr)\leq\left(\frac{a_{1}(t)M_{2}}{M_{1}}\right)^{a_{2}(t)(p+q)}M_{g_{f}}(r)
\end{equation}
where
\begin{equation}\label{a1a2}
a_{1}(t)\leq\frac{300(\sqrt{t}+1)t^{3/2}}{(t-1)^{2}},\ \ \
a_{2}(t)\leq \frac{36(\sqrt{t}+1)^{2}+162\ln\left(\frac{108e}{\sqrt{t}-1}\right)}{(\sqrt{t}-1)^{4}}.
\end{equation}
\end{Th}
\begin{R}\label{re1}
{\rm (A) As we will see from the proof inequality (\ref{cond1}) guarantees the fulfillment of conditions (1), (2) of the statement formulated in section 1.2.\\
(B) A similar to Theorem \ref{te1} result is valid for
$f$ satisfying $R_{f}(r,t,t^{2})<t$.
In this case the function $\tilde f:=f-f(0)$ satisfies (\ref{e12}) with $M_{1}:=M_{\tilde f}(r/t)$ and $M_{2}:=M_{\tilde f}(tr)$. Thus if $g$ is such that
$\tilde g\in {\cal F}_{p,q}(r;t;M_{\tilde f}(tr))$ with
$p\leq \ln\left(\frac{1+t}{2\sqrt{t}}\right)N_{\tilde f}(r,t)$, where
$\tilde g(z,w):=g(z,w+f(0))$, $(z,w)\in\Bo_{tr}^{n}\times\Di_{3M_{\tilde f}(tr)}$, then
(since $\tilde g_{\tilde f}=g_{f}$)
\begin{equation}\label{e15'}
\sup_{\Bo_{r}^{n}\times\Di_{M_{\tilde f}(tr)}}|\tilde g|\leq
\left(\frac{a_{1}(t)M_{\tilde f}(tr)}{M_{\tilde f}(r/t)}\right)^{a_{2}(t)(p+q)}M_{g_{f}}(r),
\end{equation}
\begin{equation}\label{e18'}
M_{g_{f}}(tr)\leq\left(\frac{a_{1}(t)M_{\tilde f}(tr)}{M_{\tilde f}(r/t)}\right)^{a_{2}(t)(p+q)}M_{g_{f}}(r).
\end{equation}
}
\end{R}

The proof of Theorem \ref{te1} is based on Cartan type inequalities for univariate holomorphic functions along with some geometric arguments. It is presented in sections 4-6.
In the next section we formulate several corollaries of Theorem \ref{te1} illustrated by some examples.
\sect{\hspace*{-1em}. Applications and Examples}
{\bf 2.1.} We set
\begin{equation}\label{e19}
c(M_{1},M_{2},t):=a_{2}(t)\cdot\ln\left(\frac{a_{1}(t)M_{2}}{M_{1}}\right).
\end{equation}
As a corollary of inequality (\ref{e18}) we obtain the following inequalities.
\begin{itemize}
\item[(1)] ({\em Bernstein type inequality})
\begin{equation}\label{e110}
\ln\left(\frac{M_{g_{f}}(ts)}{M_{g_{f}}(s)}\right)\leq c(M_{1},M_{2},t)(p+q), \ \ \  0<s\leq r.
\end{equation}
\item[(2)] ({\em Markov type inequality})

There is a constant $c_{1}(t)>0$ such that
\begin{equation}\label{e111}
\begin{array}{c}
\displaystyle
M_{D_{v}(g_{f})}(s)\leq\frac{c_{1}(t) \ \!c(M_{1},M_{2},t)(p+q)}{s}M_{g_{f}}(s),\\
\\
\displaystyle
\ 0<s\leq r, \ v\in\Co^{n},\ ||v||=1.
\end{array}
\end{equation}
Here $||\cdot||$ is the $l_{2}$-norm on $\Co^{n}$ and $D_{v}$
is the derivative in the direction $v$.
\item[(3)] ({\em Remez type inequality})

Consider the function $\Phi(x):=x+\sqrt{x^{2}-1}$, $|x|\geq 1$.
Then there is a constant $c_{2}(t)>0$  such that
\begin{equation}\label{e112}
\begin{array}{c}
\displaystyle
\ln M_{g_{f}}(s;z)\leq c_{2}(t)\ \!c(M_{1},M_{2},t)(p+q)\ln
\left(\Phi\left(\frac{1+\sqrt[2n]{1-\lambda}}{1-\sqrt[2n]{1-\lambda}}\right)\right)+\\
\\
\displaystyle
\sup_{\omega}\ln |g_{f}|\leq
c_{2}(t)\ \! c(M_{1},M_{2},t)(p+q)\ln\left(\frac{8n}{\lambda}\right)+
\sup_{\omega}\ln |g_{f}|
\end{array}
\end{equation}
for every Lebesgue measurable $\omega\subset\Bo_{s}^{n}(z)$ with
$\lambda:=\frac{\lambda_{2n}(\omega)}{\lambda_{2n}(\Bo_{s}^{n}(z))}$ and every ball
$\Bo_{s}^{n}(z)\subset\Bo_{r}^{n}$ (here $\lambda_{2n}$ is the Lebesgue measure on $\Co^{n}$).
\item[(4)] ({\em Jensen type inequality}) 

If $n=1$, then
\begin{equation}\label{e113}
n_{g_{f}}(r)\leq\frac{c(M_{1},M_{2},t)(p+q)}{\ln\left(\frac{1+t^{2}}{2t}\right)}.
\end{equation}
\end{itemize}
\begin{R}\label{re2}
{\rm It is well known how to derive inequalities (\ref{e110})-(\ref{e113}) from (\ref{e18}). For instance,
(\ref{e110}) and (\ref{e111}) are obtained 
by means of the Hadamard three circle inequality, see (\ref{e22}), and the Cauchy integral formula for holomorphic functions (see section 3). Inequality (\ref{e113}) is obtained by the Jensen type inequality for the number of zeros of a holomorphic function proved in [VP] (see  (\ref{vp})).
Finally, to get inequality (\ref{e112}) one repeats literally the arguments of the proof of Theorem 1.2 of [Br5] replacing inequalities (2.31) and (2.26) of [Br5] by their sharp forms presented in [BG], see there Lemmas 3 and 1.}
\end{R}
\begin{E}\label{ex1}
{\rm 
Assume that $f$ is a holomorphic homogeneous polynomial on $\Co^{n}$ of degree $d\geq 1$. Then $f$ clearly satisfies conditions (\ref{e12}) for each $r$ with $t=9$ and $M_{1}:=M_{f}(r/9)$,
$M_{2}:=M_{f}(9r)$. Now, according to (\ref{e19}), (\ref{a1a2})
and (\ref{e11}) we have for some $c_{1}<87$, $c_{2}<510$,
\begin{equation}\label{eq21}
c(M_{1},M_{2},9)< c_{1}\ln (c_{2}\cdot (81)^{d})\ \ \
{\rm and}\ \ \ N_{f}(r,9)=d.
\end{equation}
Then for a function $g\in {\cal F}_{p,q}(r;9;M_{2})$ with
$p\leq\ln(5/3)\cdot d$ Theorem \ref{te1} implies
\begin{equation}\label{eq22}
\sup_{\Bo_{r}^{n}\times\Di_{M_{2}}}|g|\leq (81)^{87(d+2)(p+q)}M_{g_{f}}(r),\ \ \ M_{g_{f}}(9r)\leq (81)^{87(d+2)(p+q)}M_{g_{f}}(r).
\end{equation}
In particular, if $g$ is a holomorphic polynomial of
degrees $k$ in $z\in\Co^{n}$ and $l$ in $w\in\Co$,
then by the classical Bernstein inequality we have
$g\in {\cal F}_{p,q}(r;9;M_{2})$ with $p:=k\ln 9$, $q:=l$.
In this case inequalities (\ref{eq22}) are valid for all $k\leq \frac{\ln(5/3)}{\ln 9}d<\frac{1}{4}d$. 

The last estimate is sharp up to an absolute factor. Indeed, for a univariate holomorphic polynomial $h$ of degree $l-1$, the polynomial $g(z,w):=(w-f(z))h(w)$ of degree $d$ in $z$ belongs to the class ${\cal F}_{p,q}(r;9;M_{2})$ with $p=d\ln 9$, $q=l$. Since $g_{f}\equiv 0$, it does not satisfy the first inequality in (\ref{eq22}). Thus inequalities (\ref{eq22}) hold for all polynomials of degrees $k$ in $z$ and $l$ in $w$, if $k<d$.

Next, for the class ${\cal F}_{d,d}(r;9;M_{2})$ we determine the
constant $\gamma_{d}(r;M_{2})$ by the formula
\begin{equation}\label{eq23}
\gamma_{d}(r;M_{2}):=\sup_{g\in {\cal F}_{d,d}(r;9;M_{2})}\left\{\sup_{\Bo_{r}^{n}\times\Di_{M_{2}}}\ln |g|-\ln M_{g_{f}}(r)\right\}
\end{equation} 
where the supremum is taken over all $g\neq 0$.\\
Observe that the polynomial $g(z,w):=w^{d}\in {\cal F}_{d,d}(r;9;M_{2})$ and satisfies
$$
\sup_{\Bo_{r}^{n}\times\Di_{M_{2}}}|g|:=(M_{f}(9r))^{d}=
9^{d^{2}}M_{g_{f}}(r).
$$
This and (\ref{eq22}) imply that
\begin{equation}\label{q24}
\ln 9\cdot d^{2}\leq\gamma_{d}(r;M_{2})\leq 348\ln 9\cdot (d+2)d.
\end{equation}

Thus in the case $p=q=d$ the logarithm of the constant in (\ref{eq22}) up to an absolute factor coincides with the
optimal constant $\gamma_{d}(r;M_{2})$.
}
\end{E}

One can also obtain analogs of inequalities (\ref{e15}),
(\ref{e18}) for
restrictions of holomorphic functions to complex submanifolds  of $\Co^{N}$ of codimension $\geq 2$. However, in general the application of Theorem \ref{te1} requires some additional conditions imposed on these submanifolds. In this paper we present
only the case of complex curves in $\Co^{N}$ for which no additional conditions are required. 

Assume that holomorphic functions $f_{i}$ on $\Di_{tr}$ satisfy conditions (\ref{e12}) with bounds
$M_{i1}$ and $M_{i2}$, $1\leq i\leq k$. We fix a permutation $\{i_{1},\dots, i_{k}\}$ of $\{1,\dots, k\}$ such that
\begin{equation}\label{e114}
N_{f_{i_{1}}}(r,t)\leq N_{f_{i_{2}}}(r,t)\leq\cdots\leq N_{f_{i_{k}}}(r,t).
\end{equation}

Let $g$ be a holomorphic function in the domain $\Di_{tr}\times\Di_{3M_{12}}\times\cdots\times\Di_{3M_{k2}}\subset\Co^{k+1}$. Suppose that for some nonnegative $p,q_{1},\dots, q_{k}$ and all $1\leq i\leq k$ 
\begin{equation}\label{e115}
\!g(\cdot,w_{1},\dots, w_{i-1},\cdot,w_{i+1},\dots, w_{k})\in {\cal F}_{p,q_{i}}(r;t;M_{i2}) \ \
{\rm for\ all} \ \ 
w_{j}\in\Di_{3M_{j2}},\ j\neq i.
\end{equation}

We set
\begin{equation}\label{e116}
\Phi(z):=(f_{1}(z),\dots, f_{k}(z))\ \ \ {\rm and}\ \ \
g_{\Phi}(z):=g(z,\Phi(z)),\ \ \ z\in\Di_{tr}.
\end{equation}

Next, we determine the sequence of nonnegative numbers $p_{0},p_{1},\dots, p_{k}$ by the formulas
\begin{equation}\label{e117}
p_{0}:=p\ \ \ {\rm and}\ \ \ p_{j}:=c(M_{j1},M_{j2},t)(p_{j-1}+q_{i_{j}}),\ \ \ 1\leq j\leq k,
\end{equation}
where $c(M_{j1},M_{j2},t)$ are defined in (\ref{e19}).
\begin{Th}\label{te12}
Assume that 
$$
p_{j}\leq\ln\left(\frac{1+t}{2\sqrt{t}}\right)\cdot N_{f_{i_{j+1}}}(r,t)\ \ \ {\rm for\ all}\ \ \ 0\leq j\leq k-1.
$$
Then
\begin{equation}\label{e118}
\max_{\Di_{r}\times\Di_{M_{12}}\times\cdots\times\Di_{M_{k2}}}|g|\leq e^{p_{1}+\cdots +p_{k}}M_{g_{\Phi}}(r)\ \ \ {\rm and}\ \ \
M_{g_{\Phi}}(tr)\leq e^{p_{k}}M_{g_{\Phi}}(r).
\end{equation}
\end{Th}
\begin{R}\label{re3}
{\rm Since $N_{f_{j}}(r):=v_{f_{j}}(r/\sqrt{t})\geq 1$ for all $1\leq j\leq k$, inequalities (\ref{e118}) are always valid for functions $g$ with sufficiently small $p,q_{0},\dots, q_{k-1}$.
}
\end{R}
{\bf 2.2.} In [CP] Coman and Poletsky obtained an important result on polynomial type inequalities for restrictions of holomorphic polynomials to certain transcendental curves in $\Co^{2}$. In this part we establish some multi-dimensional generalizations of their result that can be considered as corollaries of Theorem \ref{te1}.

Let us recall that an entire function $f$ on $\Co^{n}$ is of {\em order} $\rho\geq 0$ if
\begin{equation}\label{e129}
\rho=\limsup_{r\to\infty}\frac{\ln m_{f}(r)}{\ln r}.
\end{equation}
If $\rho<\infty$, then $f$ is called of {\em finite order}.

The following result was proved in [CP, Theorem 1.1]:\\
\\
{\bf Theorem.} {\em For any entire  function $f$ on $\Co$ of finite order $\rho>0$, there exist sequences  $\{n_{j}\}\subset\N$ convergent to $\infty$ and $\{\epsilon_{j}\}\subset\Re_{+}$ convergent to $0$ such that
for every holomorphic polynomial $g$ on $\Co^{2}$ of degree $n_{j}$ one has
\begin{equation}\label{e130}
\sup_{\Di\times\Di}|g|\leq e^{C_{1}n_{j}^{2}\ln n_{j}} M_{g_{f}}(1),
\ \ \ M_{g_{f}}(r)\leq e^{C_{2}n_{j}^{2}\ln r} M_{g_{f}}(1),\ \ \
1\leq r\leq\frac{1}{2}n_{j}^{1/\rho-\epsilon_{j}}.
\end{equation}

For every $r\geq 1$ there exists an integer $j_{r}$ such that if
$j\geq j_{r}$, then
\begin{equation}\label{e131}
n_{g_{f}}(r)\leq C_{3}n_{j}^{2},\ \ \ \frac{M_{g_{f}}(2r)}{M_{g_{f}}(r)}\leq 2^{an_{j}^{2}},\ \ \
M_{(g_{f})'}(r)\leq C_{4}n_{j}^{2}\frac{M_{g_{f}}(r)}{r}.
\end{equation}

Moreover, all the constants are effectively computed and depend only on $\rho$.}\\

The proof of this theorem is based on the Ahlfors theory of coverings surfaces and certain results of Dufresnoy along with Cartan type estimates.

Let us present a multi-dimensional generalization of this result.
\begin{Th}\label{te13}
Let $f$ be a nonpolynomial entire function on $\Co^{n}$ of order $\rho$. Then there exist sequences $\{n_{j}\}, \{r_{j}\}\subset\Re_{+}$ convergent to $\infty$ and
$\{\epsilon_{j}\}\subset\Re_{+}$ convergent to $0$
such that for every function $g\in{\cal F}_{p,q}(er_{j};e; M_{f}(e^{2}r_{j}))$ with $p\leq n_{j}$
and every $1\leq r\leq r_{j}$ the following inequalities hold:
\begin{itemize}
\item[(a)]\ \ \
$\displaystyle \sup_{\Bo^{n}\times\Di}|g|\leq e^{C_{\rho} n_{j}^{1+\epsilon_{j}}\ln r_{j}\max\{p,q\}}M_{g_{f}}(1);$
\item[(b)]\ \ \
$\displaystyle M_{g_{f}}(r)\leq e^{C_{\rho}n_{j}^{1+\epsilon_{j}}\ln r\max\{p,q\}} M_{g_{f}}(1);$
\item[(c)]\ \ \
$\displaystyle \frac{M_{g_{f}}(er)}{M_{g_{f}}(r)}\leq e^{C_{\rho}n_{j}^{1+\epsilon_{j}}\max\{p,q\}};$
\item[(d)]\ \ \
$\displaystyle
M_{D_{v}(g_{f})}(r)\leq c_{1}C_{\rho}n_{j}^{1+\epsilon_{j}}\max\{p,q\}\frac{M_{g_{f}}(r)}{r},
\ \ \ v\in\Co^{n},\ ||v||=1;$
\item[(e)]\ \ \
$\displaystyle
\ln M_{g_{f}}(s;z)\leq
c_{2}C_{\rho}n_{j}^{1+\epsilon_{j}}\max\{p,q\}\ln\left(\frac{8\lambda_{2n}(\Bo_{s}^{n})}{\lambda_{2n}(\omega)}\right)+
\sup_{\omega}\ln |g_{f}|$\ \ \
for every Lebesgue measurable set $\omega\subset\Bo_{s}^{n}(z)$ 
and every ball $\Bo_{s}^{n}(z)\subset\Bo_{r}^{n};$
\item[(f)]\ \ \
$\displaystyle n_{g_{f}}(r)\leq c_{3}C_{\rho}n_{j}^{1+\epsilon_{j}}\max\{p,q\}\ \ \ for\ \ \ n=1$.
\end{itemize}
Here $c_{1}< 9$, $c_{2}$ and $c_{3}<5$ are absolute constants and $C_{\rho}$ depends on $\rho$ only.

Moreover,
\begin{itemize}
\item[(1)]
If $\rho<\infty$, then all $\epsilon_{j}=0$ and $r_{j}\geq n_{j}^{1/(\rho+\epsilon_{j}')}$, $j\in\N$, for some sequence $\{\epsilon_{j}'\}\subset\Re_{+}$ convergent to $0$.
Also, $C_{\rho}\leq c(\ln(\rho+1)+1)^{2}(\rho+1)^{7}$ for an  absolute constant $c>0$.
\item[(2)]
If $0<\rho<\infty$, then $r_{j}\leq c_{\rho}n_{j}^{1/(\rho-\tilde\epsilon_{j}')}$,
$j\in\N$, for a sequence $\{\tilde\epsilon_{j}'\}\subset\Re_{+}$ convergent to $0$. Here $c_{\rho}\leq\left(\frac{\tilde c}{\rho_{*}}\right)^{1/\rho_{*}}$ for an absolute constant $\tilde c>0$, and $\rho_{*}:=\min\{1,\rho\}$.
\item[(3)]
If $\rho=\infty$, then  $r_{j}:=\frac{1}{e^{2}}m_{f}^{-1}(n_{j}^{1+\epsilon_{j}''})$,
$j\in\N$, for some sequence $\{\epsilon_{j}''\}\subset\Re_{+}$ convergent to $0$, and  $C_{\infty}=1$.
\end{itemize}
\end{Th}
\begin{R}\label{nr}
{\rm In the case $\rho=\infty$ we prove that $\ln r_{j}\leq n_{j}^{\delta_{j}}$ for some $\{\delta_{j}\}\subset\Re_{+}$ convergent to $0$, see (\ref{tilde}), (\ref{nj}).
Thus one can replace $n_{j}^{1+\epsilon_{j}}\ln r_{j}$ in inequality (a) by $n_{j}^{1+\tilde\epsilon_{j}}$ for some
$\{\tilde\epsilon_{j}\}\subset\Re_{+}$ convergent to $0$.
}
\end{R}
\begin{E}\label{ex2}
{\rm (A) If $g$ is a holomorphic polynomial of degree $\leq n_{j}$ on $\Co^{n+1}$, then by the classical Bernstein inequality
$g\in {\cal F}_{p,q}(er_{j};e;M_{f}(e^{2}r_{j}))$ with $p=q\leq n_{j}$. Thus Theorem \ref{te13} can be applied to such $g$.\\
(B) Let $f$ be an entire function on $\Co^{n}$ of order $1<\rho<\infty$ and $g$ be an exponential polynomial on $\Co^{n+1}$, that is,
$$
g(z,w)=\sum_{j=1}^{m}p_{j}(z,w)e^{l_{j}(z,w)},\ \ \  (z,w)\in\Co^{n}\times\Co,
$$
where $p_{j}$ is a holomorphic polynomial on $\Co^{n+1}$ of degree $d_{j}$ and $l_{j}$ is a complex linear functional on $\Co^{n+1}$ of $l_{2}$-norm $v_{j}$, $1\leq j\leq m$.

The expression
$$
m(g):=\sum_{j=1}^{m}(1+d_{j})
$$
is called the {\em degree} of $g$. Also, the {\em exponential type} of $g$ is defined by the formula
$$
\epsilon(g):=\max_{1\leq j\leq m} v_{j}.
$$

Next, let $l\in {\cal L}_{n}$ be a complex line passing through the origin. We naturally identify it with $\Co$ and define the exponential polynomial $g_{l}$ on $\Co\times\Co$ by the formula
\begin{equation}\label{e132}
g_{l}(z,w):=g(z,w),\ \ \ z\in l.
\end{equation}
Then by [VP, page 27, formula (21)] with $S^{*}:=e^{2}r_{j}$,
$S:=er_{j}$ we obtain
\begin{equation}\label{e133}
M_{g_{l}(\cdot,w)}(e^{2}r_{j})\leq e^{m(g)+2e^{2} \epsilon(g)r_{j}}M_{g_{l}(\cdot,w)}(er_{j})\ \ \ {\rm for\ all}\ \ \
w\in\Co.
\end{equation}

Similarly, by the same formula we have
\begin{equation}\label{e134}
b_{g_{l}(z,\cdot)}(3M_{f}(e^{2}r_{j}))\leq m(g)+6\epsilon(g)M_{f}(e^{2}r_{j})\ \ \ {\rm for\ all}\ \ \
z\in l.
\end{equation}

Therefore, $g\in {\cal F}_{p,q}(er_{j};e;M_{f}(e^{2}r_{j}))$ with $p:=m(g)+2e^{2}\epsilon(g)r_{j}\leq m(g)+\penalty-10000 2e^{2}c_{\rho}\epsilon(g)n_{j}^{1/(\rho-\tilde\epsilon_{j}')}$ and $q:=m(g)+6\epsilon(g)M_{f}(e^{2}r_{j})$, see Definition \ref{d11}.
Since $\rho>1$, for all sufficiently large $j$ we have $p\leq n_{j}$.
Hence for such $j$ we can apply Theorem \ref{te13}. Also, observe that $2e^{2}r_{j}\leq 6M_{f}(e^{2}r_{j})$ for all sufficiently large $j$. In particular, $\max\{p,q\}=q$ for such $j$ and inequalities of Theorem \ref{te13} are valid with $\max\{p,q\}$ substituted for
$m(g)+6\epsilon(g)M_{f}(e^{2}c_{\rho}n_{j}^{1/(\rho-\tilde\epsilon_{j}')})\leq
m(g)+\epsilon(g)e^{n_{j}^{1+\tilde\epsilon_{j}}}$ for
some $\{\tilde\epsilon_{j}\}\subset\Re_{+}$ convergent to $0$.

Suppose now that the functionals $l_{j}$ in the definition of $g$ do not depend on $w$. Thus, for a fixed $z\in\Co^{n}$, the function $g(z,\cdot)$ is a polynomial of degree $\leq d:=\max_{1\leq j\leq m}d_{j}$. In particular, instead of (\ref{e134}) we have in this case
\begin{equation}\label{e135}
b_{g_{l}(z,\cdot)}(3M_{f}(e^{2}r_{j}))\leq d.
\end{equation}
Therefore for all sufficiently large $j$ inequalities of Theorem \ref{te13} are valid with $\max\{p,q\}$ substituted for $m(g)+2e^{2}c_{\rho}\epsilon(g)n_{j}^{1/(\rho-\tilde\epsilon_{j}')}$.
}
\end{E}

Now, let us formulate some conditions under which inequalities of Theorem \ref{te13} are valid for all sufficiently large $n_{j}$ and $r_{j}$.

For a nonconstant entire function $f$ on $\Co^{n}$ of order $\rho$ we set
$$
\phi_{f}(t):=m_{f}(e^{t}),\ \ \ t\in\Re.
$$
Then $\phi_{f}$ is a convex increasing function, and so the derivative $\phi_{f}'$ exists and is continuous outside a countable set $S\subset\Re$. Also, $\phi_{f}'$ is a positive nondecreasing function on $\Re\setminus S$ having singularities of the first kind at the points of $S$. We extend $\phi_{f}'$ to $S$ by the
formula
$$
\phi_{f}'(s):=\frac{\phi_{f}'(s+)+\phi_{f}'(s-)}{2},\ \ \ s\in S,
$$
and call the extended function the derivative of $\phi_{f}$ on $\Re$.
\begin{Th}\label{te14}
Assume that $f$ satisfies one of the following conditions
\begin{itemize}
\item[(I)]
If $\rho<\infty$, 
$$
\limsup_{t\to\infty}\frac{m_{f}(e^{\alpha_{\rho}}r)-m_{f}(e^{-\alpha_{\rho}}r)+\rho e^{\rho t}}{m_{f}(e^{-\alpha_{\rho}}r)-m_{f}(e^{-2\alpha_{\rho}}r)}<A<\infty
$$
where $\alpha_{\rho}:=\min\{1,\ln(1+1/\rho)\}$.
\item[(II)]
If $\rho=\infty$, 
$$
\lim_{t\to\infty}t^{2}\left(\frac{1}{\ln \phi_{f}(t)}\right)'=0.
$$
\end{itemize}
Then there exist numbers $k_{0},r_{0}\geq 1$,
a continuous increasing to $\infty$ function 
$r:[k_{0},\infty)\to [r_{0},\infty)$ and a continuous function $\epsilon: [k_{0},\infty)\to\Re_{+}$ decreasing to $0$ as $k\to\infty$ such that
for each $k\geq k_{0}$, $r(k)\geq r_{0}$,
every 
$g\in {\cal F}_{p,q}(er(k);e;M_{f}(e^{2}r(k)))$ with $p\leq k$
and every $1\leq r\leq r(k)$ the following inequalities hold:
\begin{itemize}
\item[(a)]\ \ \
$\displaystyle \sup_{\Bo^{n}\times\Di}|g|\leq e^{C k^{1+\epsilon(k)}\ln r(k)\max\{p,q\}}M_{g_{f}}(1);$
\item[(b)]\ \ \
$\displaystyle M_{g_{f}}(r)\leq e^{Ck^{1+\epsilon(k)}\ln r\max\{p,q\}} M_{g_{f}}(1);$
\item[(c)]\ \ \
$\displaystyle \frac{M_{g_{f}}(er)}{M_{g_{f}}(r)}\leq e^{Ck^{1+\epsilon(k)}\max\{p,q\}};$
\item[(d)]\ \ \
$\displaystyle
M_{D_{v}(g_{f})}(r)\leq c_{1}Ck^{1+\epsilon(k)}\max\{p,q\}\frac{M_{g_{f}}(r)}{r},
\ \ \ v\in\Co^{n},\ ||v||=1;$
\item[(e)]\ \ \
$\displaystyle
\ln M_{g_{f}}(s;z)\leq
c_{2}Ck^{1+\epsilon(k)}\max\{p,q\}\ln\left(\frac{8\lambda_{2n}(\Bo_{s}^{n})}{\lambda_{2n}(\omega)}\right)+
\sup_{\omega}\ln |g_{f}|$\ \ \
for every Lebesgue measurable set $\omega\subset\Bo_{s}^{n}(z)$ 
and every ball $\Bo_{s}^{n}(z)\subset\Bo_{r}^{n};$
\item[(f)]\ \ \
$\displaystyle n_{g_{f}}(r)\leq c_{3}Ck^{1+\epsilon(k)}\max\{p,q\}\ \ \  for\ \ \ n=1.$
\end{itemize}
Here $c_{1}< 9$, $c_{2}$ and $c_{3}<5$ are absolute constants, for $\rho<\infty$ the constant $C$ depends on $A,\rho$ only and
$C=1$ for $\rho=\infty$.

Moreover,
\begin{itemize}
\item[(1)]
If $\rho<\infty$, then  $\epsilon\equiv 0$ and $r(k)\geq k^{1/(\rho+\epsilon'(k))}$, $k\geq k_{0}$, for some continuous function $\epsilon':[k_{0},\infty)\to\Re_{+}$ decreasing to $0$ as
$k\to\infty$.
\item[(2)]
If $0<\rho<\infty$, then $r(k)\leq ck^{1/\rho}$,
$k\geq k_{0}$,
for some $c$ depending on $A, \rho$.
\item[(3)]
If $\rho=\infty$, then  $r(k)=\frac{1}{e^{2}}m_{f}^{-1}(k^{1+\epsilon''(k)})$,
$k\geq k_{0}$, for some continuous function $\epsilon'': [k_{0},\infty)\to\Re_{+}$ decreasing to $0$ as $k\to\infty$.
\end{itemize}
\end{Th}
\begin{R}\label{rem5}
{\rm 
(A) In the case $\rho=\infty$ we show that $\ln r(k)\leq k^{\delta(k)}$ for a continuous function $\delta: [k_{0},\infty)\to\Re_{+}$ decreasing to $0$ as $k\to\infty$,
see (\ref{e751}), (\ref{e751'}).
Thus one can replace $k^{1+\epsilon(k)}\ln r(k)$ in inequality (a) by $k^{1+\tilde\epsilon(k)}$ for some continuous function
$\tilde\epsilon: [k_{0},\infty)\to\Re_{+}$ decreasing to $0$ as $k\to\infty$.\\
(B)
As an example of function $f$ satisfying condition (I) 
one can take, e.g.,
$$
f(z)=\sum_{j=1}^{m}p_{j}(z)e^{q_{j}(z)}
$$
where $p_{j},q_{j}$ are holomorphic polynomials on $\Co^{n}$.
(In this case $\lim_{r\to\infty}\frac{m_{f}(r)}{r^{\rho}}=a>0$.) 
\\
(C) As an example of function $f$ satisfying condition (II) one can take, e.g., 
$$
f(z)=e^{h(z)}\ \ \ {\rm where}\ \ \ h(z)=\sum_{j=1}^{m}p_{j}(z)e^{q_{j}(z)} 
$$
and $p_{j},q_{j}$ are holomorphic polynomials with nonnegative coefficients on $\Co^{n}$. 
}
\end{R}

Following [CP] for an entire function $f$ on $\Co^{n}$ we define
\begin{equation}\label{e136}
m_{k}(r,f):=\sup\{\ln M_{g_{f}}(r)\ :\ g\in {\cal P}_{k,n+1},\ M_{g_{f}}(1)\leq 1\},\ \ \ r\geq 1.
\end{equation}
where ${\cal P}_{k,n+1}$ is the space of holomorphic polynomials of degree $k$ on $\Co^{n+1}$.

Next, we introduce the {\em lower order of transcendence} of $f$ as
\begin{equation}\label{e137}
\underline\tau(f)=\sup\left\{\tau\ :\ \liminf_{k\to\infty}\frac{m_{k}(e,f)}{k^{\tau}}>0\right\},
\end{equation}
and
the {\em upper order of transcendence} of $f$ as
\begin{equation}\label{e138}
\overline\tau(f)=\inf\left\{\tau\ :\ \limsup_{k\to\infty}\frac{m_{k}(e,f)}{k^{\tau}}<\infty\right\}.
\end{equation}

If $f$ is a polynomial, then using the Bernstein inequality one can show that $\underline\tau(f)=\overline\tau(f)=1$.
In the case $n=1$ and $f$ is an entire function of finite positive order it was proved in [CP] that $\underline\tau(f)=2$.  Also, for each $\tau\in [3,\infty]$ there were constructed some examples of entire functions $f$ of finite positive order for which 
$\tau-1\leq\overline\tau(f)\leq\tau$.

Now, as a corollary of Theorem \ref{te13} we obtain the following generalization of the above cited result of [CP].
\begin{C}\label{c14}
If $f$ is a nonpolynomial entire function on $\Co^{n}$, then
$$
1+\frac{1}{n}\leq\underline\tau(f)\leq 2.
$$
\end{C}
\begin{R}\label{re5}
{\rm (1) Let us consider the function $\underline\tau:{\cal E}_{n}\to [1+1/n,2]$, $f\mapsto\underline\tau(f)$, defined on the set of all nonpolynomial entire functions on $\Co^{n}$. 
Since for $n=1$ the lower order of transcendence of any nonpolynomial function is $2$, one can easily construct entire functions $f$ on $\Co^{n}$, $n>1$, for which $\underline\tau(f)=2$. Thus $2$ belongs to the image of $\underline\tau$. However, we do not know what other numbers from $[1+1/n,2]$ belong to this image.\\
(2) If $f$ satisfies conditions of Theorem \ref{te14}, then
$\underline\tau(f)=\overline\tau(f)$.
}
\end{R}

In the next section we gather some auxiliary results used in the proof of Theorem \ref{te1}. Sections 4-9 are devoted to  proofs of our results.
\sect{\hspace*{-1em}. Auxiliary Results}
{\bf 3.1.} In our proofs we use the corollary of the classical Hadamard three circle inequality stating that
for a holomorphic function $h$ defined on
$\Bo_{r_{2}}^{n}$, $r_{2}>0$,
\begin{equation}\label{e22}
M_{h}(r_{1})\leq (M_{h}(r_{0}))^{1-\theta}(M_{h}(r_{2}))^{\theta},\ \ \
r_{0}\leq r_{1}\leq r_{2},\ \ \ \theta:=\frac{\ln(r_{1}/r_{0})}{\ln(r_{2}/r_{0})}.
\end{equation}
This shows that if $h\not\equiv 0$, then the function
$$
\phi_{h}(t):=m_{h}(e^{t}),\ \ \  -\infty<t<\ln r_{2},
$$
is convex and nondecreasing. In turn, the latter implies the following inequalities
\begin{itemize}
\item[(a)] For each $0<r<r_{2}$,
$$
\frac{M_{h}(r)}{M_{h}(r/e)}\leq
\frac{M_{h}(r_{2})}{M_{h}(r_{2}/e)}.
$$
\item[(b)] If $r_{2}>e$, then for each $1\leq r\leq r_{2}/e$,
$$
\frac{M_{h}(r)}{M_{h}(1)}\leq\left(\frac{M_{h}(er)}{M_{h}(r)}\right)^{\ln r}.
$$
\item[(c)] If $1<t\leq e$, then
$$
\frac{M_{h}(r_{2})}{M_{h}(r_{2}/e)}\leq\left(\frac{M_{h}(r_{2})}{M_{h}(r_{2}/t)}\right)^{\frac{1}{\ln t}}.
$$
\end{itemize}
{\bf 3.2.}
We also use Cartan type inequalities for univariate holomorphic functions. 

Let $f$ be a nonzero holomorphic function in the disk 
$\Di_{R}$. 
Fix positive $\alpha,\beta$ such that $\alpha<\beta<1$.
\begin{Th}\label{cart1}
Let $H$ be a positive number $\leq \beta e$ and 
$d>0$. Then there is a family of open disks $\{D_{j}\}_{1\leq j\leq k}$, $k\leq n_{f}(\beta R)$, 
with $\sum r_{j}^{d}\leq\frac{(2HR)^{d}}{d}$ where $r_{j}$ is the radius of 
$D_{j}$ such that
\begin{equation}\label{ca1}
\begin{array}{c}
\displaystyle
|f(z)|\geq M_{f}(\beta R)
\left(\frac{M_{f}(\alpha R)}{M_{f}(\beta R)}\right)
^{\left(\frac{\beta+\alpha}{\beta-\alpha}\right)^{2}}
\cdot\
\left(\frac{H}{\beta e}\right)^{n_{f}(\beta R)}\geq\\
\\
\\
\displaystyle 
M_{f}(\beta R)
\left(\frac{M_{f}(\alpha R)}{M_{f}(\beta R)}\right)^{\left(
\frac{\beta+\alpha}{\beta-\alpha}\right)^{2}}\cdot\ \left(
\frac{M_{f}(\beta R)}{M_{f}(R)}\right)^{
\frac{\ln\left(\frac{\beta e}{H}\right)}{\ln\left(
\frac{1+\beta^{2}}{2\beta}\right)}}\geq\\ 
\\
\\
\displaystyle
M_{f}(\beta R)
\left(\frac{M_{f}(\alpha R)}{M_{f}(R)}\right)^{\left(
\frac{\ln\left(\frac{\beta}{\alpha}\right)}{\ln
\left(\frac{1}{\alpha}\right)}\right)\cdot \left(
\frac{\beta+\alpha}{\beta-\alpha}\right)^{2}+
\frac{\ln\left(\frac{\beta e}{H}\right)}{\ln\left(
\frac{1+\beta^{2}}{2\beta}\right)}}
\end{array}
\end{equation}
for any $z\in\Di_{\alpha R}\setminus\cup_{j}D_{j}$.
\end{Th}
{\bf Proof.} We first prove the theorem for $g(z):=f(\beta Rz)$, 
$z\in\Di_{\delta}$, 
and the disks $\Di_{\gamma}\subset\Di\subset\Di_{\delta}$ where $\gamma:=\alpha/\beta$, $\delta:=1/\beta$. 

For $z\in\Di$ we write
$$
g(z):=B(z)\cdot h(z)
$$
where $B$ is the Blaschke product whose zeros are the same as for $g$ (counted with their multiplicities) and $h$ has no zeros in $\Di$. Let
$$
\rho(z,w):=\left|\frac{z-w}{1-\overline{w}z}\right|
$$
be the pseudohyperbolic metric in $\Di$. Applying to $\rho$ and $\log |B|$
the abstract Cartan estimates established in [Br6, Theorem 2.3] we have

{\em Given $H>0,\ d>0$ there is a family of open $\rho$-balls 
$\{B_{j}\}_{1\leq j\leq k}$, $k\leq n_{g}(1)=n_{f}(\beta R)$, with
$\sum r_{j}^{d}\leq\frac{(2H)^{d}}{d}$
where $r_{j}$ is the radius of $B_{j}$ such that for any 
$z\in\Di\setminus\cup_{j}B_{j}$}
\begin{equation}\label{ca2}
|B(z)|\geq \left(\frac{H}{e}\right)^{n_{f}(\beta R)} .
\end{equation}
Since each $B_{i}$ is the subset of the Euclidean
disk $D_{i}$ centered at the same point and of the same radius, the above
inequality is also valid for each $z\in\Di\setminus\cup_{j}D_{j}$. 

Next, we have $M_{f}(\beta R)=M_{h}(1)$ and
$M_{f}(\alpha R)=M_{g}(\gamma)$. These identities imply that
$M_{h}(\gamma)\geq M_{f}(\alpha R)$ and that the function $u(z):=-\ln |h(z)|+\ln M_{f}(\beta R)$ is nonnegative harmonic in $\Di$. We will apply to $u$ the classical Harnack inequality.

Take $w=\gamma e^{i\phi}$ such that $M_{h}(\gamma)=|h(w)|$ and
let
$$
G(z)=\frac{z+w}{1+\overline{w}z}
$$
be the M\"{o}bius transformation of $\Di$ sending $0$ to $w$. Then
$u_{G}(z):=u(G(z))$ is a nonnegative harmonic function in $\Di$ and $u_{G}(0)\leq
\ln[M_{f}(\beta R)/M_{f}(\alpha R)]$. By the Harnack inequality we have
$$
u(0)=u_{G}(-w)\leq u_{G}(0)\frac{1+|w|}{1-|w|}\leq
\left[\ln\left(\frac{M_{f}(\beta R)}{M_{f}(\alpha R)}\right)\right]\left(\frac{1+\gamma}{1-\gamma}\right) .
$$
Applying again the Harnack inequality to $u$ at the points 0 and
$y$ such that $|y|=\gamma$ and $u(y)=\sup_{\Di_{\gamma}}u$ and using
the previous estimate we have
$$
\sup_{\Di_{\gamma}}u\leq\left[\ln\left(\frac{M_{f}(\beta R)}{M_{f}(\alpha R)}\right)\right]
\left(\frac{1+\gamma}{1-\gamma}\right)^{2} .
$$
From here and the definition of $u$ it follows that for any $z\in\Di_{\gamma}$
\begin{equation}\label{ca3}
|h(z)|\geq M_{f}(\beta R)\cdot\left(
\frac{M_{f}(\alpha R)}{M_{f}(\beta R)}\right)^{\left(
\frac{\beta+\alpha}{\beta-\alpha}\right)^{2}} .
\end{equation}
Combining inequalities (\ref{ca2}), (\ref{ca3}) and going back to $f$ we
obtain the first inequality of (\ref{ca1}). To obtain the second inequality we use the estimate from [VP, Lemma 1]:
\begin{equation}\label{vp}
n_{f}(\beta R)\leq\frac{\ln\left(\frac{M_{f}(R)}{M_{f}(\beta R)}\right)}{
\ln\left(\frac{1+\beta^{2}}{2\beta}\right)}\ .
\end{equation}
Finally, the third inequality is obtained by the application of
the Hadamard three circle inequality estimating $M_{f}(\beta R)$ by $M_{f}(R)$ and $M_{f}(\alpha R)$ to the second term of the second line of (\ref{ca1}) and by replacing $M_{f}(\beta R)$ to $M_{f}(\alpha R)$ in the third term of that line.
\ \ \ \ \ $\Box$

Applying Theorem \ref{cart1} to $R:=tr$, $\beta R:=\sqrt{t}r$ and $\alpha R:= r$,
$r>0$, $1<t\leq 9$, we obtain
\begin{Th}\label{cart2}
Let $f$ be a nonzero holomorphic function in $\Di_{tr}$, $r>0$, $1\leq t\leq 9$.
Let $H$ be a positive number $\leq e/\sqrt{t}$. Then there is a 
family of open disks $\{D_{j}\}_{1\leq j\leq k}$, $k\leq n_{f}(\sqrt{t}r)$, with $\sum r_{j}\leq 2Htr$ where $r_{j}$ is the radius of $D_{j}$ such that for each $z\in\Di_{r}\setminus \cup_{j}D_{j}$
\begin{equation}\label{ca4}
|f(z)|\geq M_{f}(\sqrt{t}r)\left(\frac{M_{f}(r)}{M_{f}(tr)}\right)^{c(H)}
\end{equation}
where
\begin{equation}\label{ca5}
c(H):=\frac{(\sqrt{t}+1)^{4}+18(\sqrt{t}+1)^{2}\ln\left(\frac{e}{H}\right)}{2(t-1)^{2}}
\end{equation}
and
\begin{equation}\label{ca6}
n_{f}(\sqrt{t}r)\leq \frac{\ln\left(\frac{M_{f}(tr)}{M_{f}(\sqrt{t}r)}\right)}{
\ln\left(\frac{1+t}{2\sqrt{t}}\right)}\leq\frac{9(\sqrt{t}+1)^{2}\ln\left(\frac{M_{f}(tr)}{M_{f}(\sqrt{t}r)}\right)}{(t-1)^{2}}.
\end{equation}
\end{Th}
{\bf Proof.} Inequality (\ref{ca4}) follows directly from (\ref{ca1}) with 
$$
\frac{1}{2}\cdot \left(
\frac{\sqrt{t}+1}{\sqrt{t}-1}\right)^{2}+
\frac{\ln\left(\frac{e}{\sqrt{t}H}\right)}{\ln\left(
\frac{1+t}{2\sqrt{t}}\right)}
$$
instead of $c(H)$ defined by (\ref{ca5}). Here
\begin{equation}\label{ca8}
\begin{array}{c}
\displaystyle
\left(
\frac{\sqrt{t}+1}{\sqrt{t}-1}\right)^{2}=\frac{(\sqrt{t}+1)^{4}}{(t-1)^{2}}\ \ \ {\rm and}\\
\\
\displaystyle
\frac{1}{\ln\left(
\frac{1+t}{2\sqrt{t}}\right)}=\frac{1}{\ln\left(1+\frac{(\sqrt{t}-1)^{2}}{2\sqrt{t}}\right)}\leq
\frac{1}{\ln\left(1+\frac{(\sqrt{t}-1)^{2}}{6}\right)}\leq\frac{9(\sqrt{t}+1)^{2}}{(t-1)^{2}}.
\end{array}
\end{equation}
We used that $\ln(1+x)\geq\frac{2}{3}x$ for $0\leq x\leq\frac{2}{3}$. 

Now, from (\ref{ca8}) we obtain (\ref{ca4}) with $c(H)$ given by (\ref{ca5}) and inequality (\ref{ca6}).\ \ \ \ \ $\Box$
\begin{C}\label{cart3}
Under the assumptions of Theorem \ref{cart2} there exists a circle
$S_{l}:=\{z\in\Co\ :\ |z|=l\}$, $r/\sqrt{t}\leq l\leq r$, such that for each $z\in S_{l}$
\begin{equation}\label{ca9}
|f(z)|\geq M_{f}(\sqrt{t}r)\left(\frac{M_{f}(r)}{M_{f}(tr)}\right)^{\gamma(t)}
\end{equation}
where
\begin{equation}\label{ca10}
\gamma(t):=\frac{(\sqrt{t}+1)^{4}+18(\sqrt{t}+1)^{2}\ln\left(\frac{4et^{3/2}}{\sqrt{t}-1}\right)}{2(t-1)^{2}}>0.
\end{equation}
\end{C}
{\bf Proof.} We apply Theorem \ref{cart2} with $H=\frac{\sqrt{t}-1}{4t^{3/2}}$. Then the sum
of radii of the disks $D_{j}$ is $\leq\frac{\sqrt{t}-1}{2\sqrt{t}}r$. In particular, the
projection of $\cup_{j}D_{j}$ onto the radial axis (in polar
coordinates of $\Co$) is an open set of linear measure $\leq \frac{\sqrt{t}-1}{\sqrt{t}}r$. Therefore this set cannot cover the closed interval $\{s\in\Re_{+}\ :\ r/\sqrt{t}\leq s\leq r\}$.
This implies that there is a circle $S_{l}$ with $r/\sqrt{t}\leq l\leq r$ which does not
intersect $\cup_{j}D_{j}$. According to (\ref{ca4})
$f|_{S_{l}}$ satisfies the required estimate.\ \ \ \ \ $\Box$\\
\\
{\bf 3.3.} In the proofs we use also the following Markov type inequality.
\begin{Th}\label{markov}
Assume that $h$ is a holomorphic function in the ball $\Bo_{tR}^{n}$, $R>0$, $1<t\leq 9$, satisfying for some $d\geq 0$
\begin{equation}\label{e31}
M_{h}(tR)\leq e^{d}M_{h}(R).
\end{equation}
Then
\begin{equation}\label{e32}
M_{D_{v}h}(R)\leq\frac{\kappa(d;t)}{R}M_{h}(R)
\end{equation}
where $D_{v}$ is the derivative in the direction $v\in\Co^{n}$,
$||v||=1$, and
\begin{equation}\label{kap}
\kappa(d;t):=\left\{
\begin{array}{ccc}
\displaystyle
\frac{e}{t^{\ln(1+1/d)}-1}&{\rm if}&
\displaystyle d\geq\frac{1}{e-1}\\
\\
\displaystyle
\frac{2d}{\sqrt{t}-1}&{\rm if}&
\displaystyle 0\leq d<\ln\left(\frac{1+t}{2\sqrt{t}}\right)\\
\\
\displaystyle
\frac{e^{d}}{t-1}&{\rm if}&
\displaystyle
\ln\left(\frac{1+t}{2\sqrt{t}}\right)\leq d<\frac{1}{e-1}.
\end{array}
\right.
\end{equation}
\end{Th}
(Observe that $\ln\left(\frac{1+9}{2\sqrt{9}}\right)<\frac{1}{e-1}$ so that formula (\ref{kap}) is correct.)\\
{\bf Proof.} Without loss of generality we may assume that $h$ is not identically zero. We will consider several cases.

(1) Assume that $d\geq\frac{1}{e-1}$. 
Take $x\in\partial B_{R}^{n}$ and let $l=\{x+zv\ :\ z\in\Co\}$ be the complex line passing through $x$. We set
$D_{s}=\Bo_{sR}^{n}\cap l$, $1\leq s\leq t$. Then $D_{s}\subset l$ is the disk of radius $r_{s}$ centered at $c\in\Bo_{R}^{n}$ where $c$ is such that $h:=dist(l,0)=||c||$ and $r_{s}:=\sqrt{(sR)^{2}-h^{2}}$. We will naturally identify $D_{s}$ with $\Di_{r_{s}}$. It is easy to check that for
all $s\geq q\geq 1$ the following inequalities hold:
\begin{equation}\label{compare}
\frac{r_{s}}{r_{q}}\geq\frac{s}{q}\ \ \ {\rm and}\ \ \
r_{s}-r_{q}\geq (s-q)R.
\end{equation}

We set $s:=t^{\ln(1+1/d)}$ so that $1<s\leq t$. Then by means of the Hadamard three circle inequality, see (\ref{e22}), we obtain
\begin{equation}\label{eq213}
M_{h}(sR)\leq e^{(d\ln s)/(\ln t)}M_{h}(R)=
e^{d\ln(1+1/d)}M_{h}(R)\leq eM_{h}(R).
\end{equation}

Consider the function $\tilde h:=h|_{D_{t}}$ and the disk $D\subset l$ centered at $x\in D_{1}$ of radius $(s-1)R$. By (\ref{compare}) we have
$$
r_{1}+(s-1)R\leq r_{s}\leq s.
$$
Thus $D$ belongs to $D_{s}\subset\Bo_{sR}^{n}$. Now, from the Cauchy integral formula for the derivative of $\tilde h$ in $D$ by (\ref{eq213}) we get
\begin{equation}\label{eq214}
|(D_{v}h)(x)|:=|\tilde h'(x)|\leq\frac{1}{(s-1)R}
M_{h}(sR)\leq
\frac{e}{(t^{\ln(1+1/d)}-1)R}M_{h}(R).
\end{equation}

(2) Suppose now that $d<\ln\left(\frac{1+t}{2\sqrt{t}}\right)$.
Let $Z_{h}$ be the zero set of $h$. We first prove
\begin{Lm}\label{iter}
Under the above condition $Z_{h}\cap\Bo_{\sqrt{t}R}^{n}=\emptyset$.
\end{Lm}
{\bf Proof.} Assume, on the contrary, that there is $y\in Z_{h}\cap\Bo_{\sqrt{t}R}^{n}$. Take $v\in\overline{\Bo}_{R}^{n}$ such that
$$
|h(v)|=\sup_{\Bo_{R}^{n}}|h|.
$$
Let $l$ be a complex line passing through $v$ and $y$. As before we set $D_{s}=\Bo_{sR}^{n}\cap l$ and identify it with $\Di_{r_{s}}$ with an appropriate definition $r_{s}$ (see case (1)). Then for the function
$\tilde h:=h|_{D_{t}}$ we have by [VP, Lemma 1]
$$
n_{\tilde h}(r_{\sqrt{t}})
\leq\frac{\ln\left(\frac{M_{\tilde h}(r_{t})}{M_{\tilde h}(r_{\sqrt{t}})}\right)}{\ln\left(\frac{1+(r_{t}/r_{\sqrt{t}})^{2}}{2(r_{t}/r_{\sqrt{t}})}\right)}\leq
\frac{\ln\left(\frac{M_{h}(tR)}{M_{ h}(R)}\right)}{\ln\left(\frac{1+(r_{t}/r_{\sqrt{t}})^{2}}{2(r_{t}/r_{\sqrt{t}})}\right)}<
\frac{\ln\left(\frac{1+t}{2\sqrt{t}}\right)}{\ln\left(\frac{1+(r_{t}/r_{\sqrt{t}})^{2}}{2(r_{t}/r_{\sqrt{t}})}\right)}\leq 1.
$$
We used here that $M_{\tilde h}(r_{\sqrt{t}})\geq M_{h}(R)$ (by the choice of $l$), the function $x\mapsto\ln\left(\frac{1+x^{2}}{2x}\right)$ is increasing for $x\geq 1$ and $\sqrt{t}\leq\frac{r_{t}}{r_{\sqrt{t}}}$, see (\ref{compare}).

Thus $\tilde h$ has no zeros in $D_{\sqrt{t}}$. This contradicts to the assumption $y\in Z_{h}\cap D_{\sqrt{t}}$.\ \ \ \ \ $\Box$

Continuing the proof of the theorem consider the line $l$ as in the proof of case (1). Then according to the lemma the corresponding function $\tilde h:=h|_{D_{t}}$ has no zeros on $\Di_{r_{\sqrt{t}}}$. In particular, the holomorphic function $g:=\ln(\tilde h/M_{h}(R))$ is well defined there (for some choice of the branch of the logarithm).  Also, the function $g+\overline{g}=\ln|\tilde h/M_{ h}(R)|^{2}$ is harmonic on $\Di_{r_{\sqrt{t}}}$. Now from the Cauchy integral formula in the disk centered at $x\in\Di_{r_{1}}$ of radius $(\sqrt{t}-1)R$ we obtain
$$
\frac{|\tilde h'(x)|}{|\tilde h(x)|}:=|g'(x)|\leq\frac{\sup_{r_{\sqrt{t}}}|g+\overline{g}|}{(\sqrt{t}-1)R}\leq\frac{\sup_{\sqrt{t}R}\ln| h/M_{h}(R)|^{2}}{(\sqrt{t}-1)R}
\leq\frac{2d}{(\sqrt{t}-1)R}.
$$
(We used here that $r_{1}+(\sqrt{t}-1)R\leq r_{\sqrt{t}}\leq \sqrt{t}R$, see (\ref{compare}).)

This implies
\begin{equation}\label{eq215}
|(D_{v}h)(x)|:=|\tilde h'(x)|\leq\frac{2d}{(\sqrt{t}-1)R}M_{h}(R).
\end{equation}

(3) Finally, assume that $\ln\left(\frac{1+t}{2\sqrt{t}}\right)\leq d<\frac{1}{e-1}$.
Applying to $\tilde h:=h|_{D_{t}}$ the Cauchy integral formula for $x\in D_{1}\subset l$ (with $l$ as in case (1)) we have:
\begin{equation}\label{eq216}
\begin{array}{c}
\displaystyle
|(D_{v}h)(x)|:=|\tilde h'(x)|\leq\frac{1}{(t-1)R}M_{\tilde h}(r_{t})\leq\frac{e^{d}}{(t-1)R}M_{h}(R)\leq\\
\\
\displaystyle
\frac{e^{1/(e-1)}d}{(t-1)\ln\left(\frac{1+t}{2\sqrt{t}}\right)R}M_{h}(R).
\end{array}
\end{equation}
Inequalities (\ref{eq214}), (\ref{eq215}), (\ref{eq216}) imply inequality (\ref{e32}).\ \ \ \ \ $\Box$
\begin{R}\label{emark}
{\rm (1) If $t=e$, then 
\begin{equation}\label{kap1}
\kappa(d;e)\leq\max\left\{ed,\frac{2d}{\sqrt{e}-1},\frac{e^{1/(e-1)}d}{(e-1)\ln\left(\frac{1+e}{2\sqrt{e}}\right)}\right\}<9d.
\end{equation}
(2) For $d\geq\frac{1}{e-1}$ we have by the mean-value inequality for $f(x):=x^{\ln t}$,
\begin{equation}\label{case1}
\kappa(d;t):=\frac{e}{t^{\ln(1+1/d)}-1}=
\frac{e}{(1+1/d)^{\ln t}-1}\leq\left\{
\begin{array}{ccc}
\displaystyle \frac{ed}{\ln t}&{\rm if}&t\geq e\\
\\
\displaystyle \frac{e^{2}d}{t\ln t}&{\rm if}&t<e.
\end{array}
\right.
\end{equation}
}
\end{R}
\sect{\hspace*{-1em}. A Geometric Result} {\bf 4.1.} The proof of Theorem \ref{te1} is based on the following result.

Let $F$ be a nonconstant holomorphic function in $\Di_{t}$, $1<t\leq 9$, satisfying
\begin{equation}\label{e21}
M_{F}(t)\leq 1,\ \ \ M_{F}(1/t)\geq M\ \ \ {\rm and}\ \ \
\frac{M_{F}(1/t)}{M_{F}(1/t^{2})}\geq\sqrt{t}.
\end{equation}

We set 
\begin{equation}\label{set}
\begin{array}{c}
\displaystyle
N_{F}(t)=N_{F}(1,t):=v_{F}(1/\sqrt{t}),\ \ \
\lambda(t):=\frac{9(\sqrt{t}+1)^{2}}{(t-1)^{2}}\ln\left(\frac{2(\sqrt{t}+1)}{M(t-1)^{2}}\right),\\
\\
\displaystyle
\gamma(t):=\frac{(\sqrt{t}+1)^{4}+18(\sqrt{t}+1)^{2}\ln\left(\frac{4et^{3/2}}{\sqrt{t}-1}\right)}{2(t-1)^{2}},\ \ \ r_{0}(t):=\left(\frac{M(t-1)}{4(\sqrt{t}+1)}\right)^{\gamma(t)+1}.
\end{array}
\end{equation}
\begin{Th}\label{te2}
There is a number $c\in\Co$, $|c|<1$, and for each $y\in\Co$,
$|y|\leq r_{0}(t)$, and $s\in (0,r_{0}(t)/3]$ there is
$c_{y,s}\in\Co$, $|c_{y,s}|<s$, such that the set of zeros of the
function $F-c-y-c_{y,s}$ in $\Di$ contains at least $N_{F}(t)$ points
with pairwise distances greater than $\frac{s(t-1)}{\sqrt{\lambda(t)}}$.
\end{Th}

This result can be reformulated as follows.\\

{\em Let $\Gamma:=\{(z,F(z))\in\Co^{2}\ :\ z\in\Di_{t}\}$ be the graph of $F$. There is a number $c\in\Co$, $|c|<1$, such that for each point $v=(x,c+y)\in\Di_{t}\times \overline{\Di}_{r_{0}(t)}(c)$
and every $s\in (0, r_{0}(t)/3]$, there is a point $v'=(x,c+y+c_{y,s})\in\{x\}\times\Co$, $||v'-v||<s$, such that the complex line $l:=\{(z,w)\in\Co^{2}\ :\ w-c-y-c_{y,s}=0\}$, parallel to the $z$-axis and passing through $v'$, intersects the graph $\Gamma$ over $\Di$ in at least $N_{F}(t)$ points with pairwise distances greater than $\frac{s(t-1)}{\sqrt{\lambda(t)}}$.}\\

In sections 4.2-4.4 we formulate some auxiliary results used in the proof of Theorem \ref{te2}.\\
\\
{\bf 4.2.} Applying the Hadamard three circle inequality (\ref{e22}) to our function $F$ with $r_{0}:=1/t^{2}$,
$r_{1}:=1/t$ and $r_{2}:=1$ from (\ref{e21}) we obtain
\begin{equation}\label{e21'}
\frac{M_{F}(1)}{M_{F}(1/t)}\geq\frac{M_{F}(1/t)}{M_{F}(1/t^{2})}\geq\sqrt{t}.
\end{equation}
Then applying (\ref{e22}) with $r_{0}:=1/t$,
$r_{1}:=1$ and $r_{2}:=t$ from (\ref{e21'}) we obtain 
\begin{equation}\label{e23}
\frac{1}{M}\geq\frac{M_{F}(t)}{M_{F}(1/t)}\geq t.
\end{equation}

We use this estimate to prove 
\begin{Lm}\label{le1}
For $F':=\frac{dF}{dz}$ we have
\begin{equation}\label{e24}
n_{F'}(1)<\lambda(t).
\end{equation}
\end{Lm}
{\bf Proof.}
For $z\in\Di$, using the mean-value inequality
$|F(z)-F(0)|\leq M_{F'}(1)$, we obtain
\begin{equation}\label{e25}
M_{F'}(1)\geq M_{F}(1)-M_{F}(1/t).
\end{equation}
Also, by the Cauchy integral formula for $F'$ we get
\begin{equation}\label{e26}
M_{F'}(\sqrt{t})\leq\sup_{z\in\Di_{\sqrt{t}}}\left\{\frac{1}{2\pi}\int_{S_{t-\sqrt{t}}(z)}\frac{|F(\xi)|}{|\xi-z|^{2}}|d\xi|\right\}
\leq\frac{M_{F}(t)}{\sqrt{t}(\sqrt{t}-1)}<\frac{2}{t-1}.
\end{equation}
Here $S_{t-\sqrt{t}}(z)$ stands for the boundary of the disk $\Di_{t-\sqrt{t}}(z)$.

Finally, we apply the Jensen inequality for the number of zeros of a
holomorphic function proved in [VP]. Then from (\ref{e25}),
(\ref{e26}), (\ref{e21}), (\ref{e21'}) and (\ref{ca8}) we obtain
$$
\begin{array}{c}
\displaystyle
n_{F'}(1)\leq\frac{1}{\ln\left(\frac{1+t}{2\sqrt{t}}\right)}
\ln\left(\frac{M_{F'}(\sqrt{t})}{M_{F'}(1)}\right)<
\frac{9(\sqrt{t}+1)^{2}}{(t-1)^{2}}\ln\left(\frac{2}{(t-1)(M_{F}(1)-M_{F}(1/t))}\right)\leq\\
\\
\displaystyle
\frac{9(\sqrt{t}+1)^{2}}{(t-1)^{2}}\ln\left(\frac{2(\sqrt{t}+1)}{M(t-1)^{2}}\right).
\ \ \ \ \ \Box
\end{array}
$$
{\bf 4.3.} By the definition of $N_{F}(t)$ there is a number $c\in\Co$, $|c|\leq M_{F}(1/\sqrt{t})$, such that the function
$F_{c}:=F-c$ has $N_{F}(t)$ zeros in $\Di_{1/\sqrt{t}}$. For this function
we have
\begin{equation}\label{e27}
M_{F_{c}}(t)\leq M_{F}(t)+M_{F}(1/\sqrt{t})<2.
\end{equation}

Further, $M_{F_{c}}(1)>M_{F_{c}}(1/t)\geq |c|-M_{F}(1/t)$.
Assuming, first, that $|c|\geq\frac{1+\sqrt{t}}{2\sqrt{t}}M_{F}(1)$ and using (\ref{e21'}) we get from here 
$$
M_{F_{c}}(1)>\frac{1+\sqrt{t}}{2\sqrt{t}}M_{F}(1)-M_{F}(1/t)\geq\frac{\sqrt{t}-1}{2}M_{F}(1/t)\geq\frac{\sqrt{t}-1}{2}M.
$$
Assume now that $|c|<\frac{1+\sqrt{t}}{2\sqrt{t}}M_{F}(1)$. Then
$$
M_{F_{c}}(1)\geq M_{F}(1)-|c|>\frac{\sqrt{t}-1}{2\sqrt{t}}M_{F}(1)\geq\frac{\sqrt{t}-1}{2}M_{F}(1/t)\geq
\frac{\sqrt{t}-1}{2}M.
$$
Thus we have
\begin{equation}\label{e28}
M_{F_{c}}(1)>\frac{\sqrt{t}-1}{2}M.
\end{equation}
In particular, from (\ref{e27}) and (\ref{e28}) we obtain
\begin{equation}\label{e29}
\frac{M_{F_{c}}(t)}{M_{F_{c}}(1)}\leq \frac{4(\sqrt{t}+1)}{M(t-1)}.
\end{equation}
From here and the Jensen inequality of [VP] we get (recall that
$1<t\leq 9$)
\begin{equation}\label{e29'}
\begin{array}{c}
\displaystyle
\!N_{F}(t)\leq n_{F_{c}}(1)\leq\frac{1}{\ln\left(\frac{1+t^{2}}{2t}\right)}\ln\left(\frac{M_{F_{c}}(t)}{M_{F_{c}}(1)}\right)\leq\\
\\
\displaystyle
\frac{9(\sqrt{t}+1)^{2}}{2(t-1)^{2}}\ln\left(\frac{4(\sqrt{t}+1)}{M(t-1)}\right)=:\delta(t).
\end{array}
\end{equation}
{\bf 4.4.} We will also use Corollary \ref{cart3}. According to this corollary for $f:=F_{c}$ and $r:=1$ using
(\ref{e28}), (\ref{e29}) we obtain that

{\em there is a circle $S_{l}$ with $1/\sqrt{t}\leq l\leq 1$ such that}
\begin{equation}\label{e213}
|F_{c}(z)|>2\left(\frac{M(t-1)}{4(\sqrt{t}+1)}\right)^{\gamma(t)+1}=:2r_{0}(t)\ \ \ {\rm for\ all}\ \ \
z\in S_{l}.
\end{equation}
{\bf 4.5. Proof of Theorem \ref{te2}.} Let $c\in\Co$, $|c|<1$, be
the number introduced in section 3.3. We will prove that for each
$y\in\Co$, $|y|\leq r_{0}(t)$, and $s\in (0, r_{0}(t)/3]$ there is
$c_{y,s}\in\Co$, $|c_{y,s}|<s$, such that the set of zeros of the
function $F-c-y-c_{y,s}$ in $\Di$ contains at least $N_{F}(t)$ points
(see (\ref{set})) with pairwise distances greater than $\frac{s(t-1)}{\sqrt{\lambda(t)}}$.

First, from inequality (\ref{e213}) by the Rouch\'{e} theorem we
deduce that
\begin{equation}\label{e214}
n_{F_{c}-a}(l)=n_{F_{c}}(l)\geq N_{F}(t) \ \ {\rm for\ all}\ \ \
a\in\Co,\ |a|\leq 2r_{0}(t).
\end{equation}
This is valid, in particular, for $a:=y+b$ with $|b|\leq r_{0}(t)$.

Let $C_{F}\subset F_{c}(\Di)\subset\Co$ be the set of critical
values of $F_{c}|_{\Di}$.
\begin{Lm}\label{le3}
For each $s\in (0,r_{0}/3]$ there is $c_{y,s}\in\Co$, $|c_{y,s}|<s$, such
that
$$
dist(y+c_{y,s},C_{F})>\frac{s}{\sqrt{\lambda(t)}}.
$$
\end{Lm}
(Observe that by (\ref{e23}) and (\ref{set}), $\sqrt{\lambda(t)}>\frac{1}{2}$ for $1<t\leq 9$.)\\
{\bf Proof.} 
We will assume that $C_{F}\neq\emptyset$. For otherwise, we set $c_{y,s}=0$.

By Lemma \ref{le1} the number of critical points of $F_{c}$
in $\Di$ is $<\lambda(t)$. Since $N_{F}(t)\geq 1$, 
(\ref{e214}) implies that $\Di_{s}(y)\subset
F_{c}(\Di)$. In particular, $\Di_{s}(y)\cap C_{F}$ contains $<\lambda(t)$ points. Thus there is $c_{y,s}\in\Di_{s}$ such that
$dist(y+c_{y,s},C_{F})>\frac{s}{\sqrt{\lambda(t)}}$. Indeed, for
otherwise, the closed disks of radius $\frac{s}{\sqrt{\lambda(t)
}}$ centered at the points of $\Di_{s}(y)\cap C_{F}$ cover
$\Di_{s}(y)$. Comparing the areas of $\Di_{s}(y)$ and of this cover we
obtain a contradiction:
$$
\pi s^{2}<\lambda(t)\cdot \pi\cdot\left(\frac{s}{\sqrt{\lambda(t)}}\right)^{2}=\pi s^{2}.\ \ \ \ \ \Box
$$

Now, by Lemma \ref{le3} we obtain that $\Di_{r_{1}}(y+c_{y,s})\cap C_{F}=\emptyset$,
$r_{1}:=\frac{s}{\sqrt{\lambda(t)}}$. Moreover, by (\ref{e214}) 
we have $\Di_{r_{1}}(y+c_{y,s})\subset F_{c}(\Di_{l})$ because
$|c_{y,s}|+r_{1}<3s\leq r_{0}(t)$. Thus from (\ref{e214}) for
$X:=F_{c}^{-1}(\Di_{r_{1}}(y+c_{y,s}))\cap\Di_{l}$ we obtain that $F_{c}:X\to\Di_{r_{1}}(y+c_{y,s})$ is a proper conformal map and
$\#\{F_{c}^{-1}(z)\cap\Di_{l}\}$ is the same for any
$z\in\Di_{r_{1}}(y+c_{y,s})$. Hence
$F_{c}:X\to\Di_{r_{1}}(y+c_{y,s})$ is an unbranched covering of
$\Di_{r_{1}}(y+c_{y,s})$ consisting of at least $N_{F}(t)$ sheets. In
particular, $X$ is biholomorphic to the disjoint union of $k$ copies
of $\Di_{r_{1}}(y+c_{y,s})$ where $k$ is the number of sheets of
$F_{c}|_{X}$.

We set $Y=\{y_{1},\dots, y_{k}\}:=F_{c}^{-1}(y+c_{y,s})\cap\Di_{l}$.
Assume that for some $i\neq j$ we have
$|y_{i}-y_{j}|\leq (t-1)r_{1}$. Then by the mean-value theorem for
each $z$ from the interval $\gamma:=[y_{i},y_{j}]$ we have
$$
|F_{c}(z)-F_{c}(y_{i})|\leq
M_{F'}(l)|y_{i}-y_{j}|<\frac{M_{F}(t)}{t-1}(t-1)r_{1}\leq r_{1}.
$$
Thus $F_{c}(z)\in \Di_{r_{1}}(y+c_{y,s})$. In particular,
$F_{c}(\gamma)\subset \Di_{r_{1}}(y+c_{y,s})$ is a closed curve.
Then since $F_{c}:X\to\Di_{r_{1}}(y+c_{y,s})$ is an unbranched
covering, $\gamma$ should be a closed curve, as well. This
contradiction shows that $|y_{i}-y_{j}|>(t-1)r_{1}$ for all
$i\neq j$.

The proof of Theorem \ref{te2} is complete.\ \ \ \ \ $\Box$
\sect{\hspace*{-1em}. Proof of Theorem \ref{te1}: Case $n=1$} 
{\bf 5.1.} First, we will prove Theorem \ref{te1} for the functions $f$ and $g$ satisfying the assumptions of the theorem for $n=1$ and such that in (\ref{e12})
\begin{equation}\label{e330}
R_{f}(r,t,t^{2}):=\frac{M_{f}(r/t)}{M_{f}(r/t^{2})}\geq\sqrt{t}.
\end{equation}

We also set
\begin{equation}\label{e33''}
M:=\frac{M_{1}}{M_{2}}.
\end{equation}

Next, we define new functions $F$ and $G$ by the
formulas
\begin{equation}\label{e30}
\begin{array}{c} \displaystyle
 F(z)=\frac{f(rz)}{M_{2}},\ \ \ z\in\Di_{t},\ \ \ {\rm and}\\
 \\
 \displaystyle
G(z,w):=g(rz, M_{2}w),\ \ \
(z,w)\in\Di_{t}\times\Di_{3}.
\end{array}
\end{equation}
Then $F$ satisfies conditions (\ref{e21}) and
$G$ satisfies the conditions 
\begin{equation}\label{e12'}
M_{G(\cdot,w)}(t)\leq e^{p}\cdot M_{G(\cdot,w)}(1)\ \ \
{\rm for\ all}\ \ \ w\in\Di_{3};
\end{equation}
\begin{equation}\label{e13'}
b_{G(z,\cdot)}(3)\leq q\ \ \ {\rm for\ all}\ \ \
z\in\Di_{t}.
\end{equation}

Now, we will prove the following version of Theorem \ref{te1}.
\begin{Th}\label{te1'}
Assume that
$p\leq\ln\left(\frac{1+t}{2\sqrt{t}}\right)\cdot N_{F}(t)$. 
Then for $G_{F}(z):=G(z,F(z))$
$$
\sup_{\Di\times\Di}|G|\leq
\left(\frac{c_{1}(t)}{M}\right)^{c_{2}(t)(p+q)}M_{G_{F}}(1)
$$
\end{Th}
where
\begin{equation}\label{c1c2}
c_{1}(t):=\frac{50(\sqrt{t}+1)}{(t-1)^{2}},\ \ \
c_{2}(t):=\frac{18(\sqrt{t}+1)^{2}+81\ln\left(\frac{108e}{\sqrt{t}-1}\right)}{(\sqrt{t}-1)^{4}}.
\end{equation}

Going back to the functions $f$ and $g$ and noticing that
$G_{F}(z)=g_{f}(rz)$, $z\in\Di_{t}$,
$\sup_{\Di\times\Di}|G|=\sup_{\Di_{r}\times\Di_{M_{2}}}|g|$, $M_{G_{F}}(1)=M_{g_{f}}(r)$ and $N_{F}(t)=N_{f}(r;t)$
we obtain from this theorem inequality (\ref{e15}) in the case $n=1$.\\
\\
{\bf 5.2. Proof of Theorem \ref{te1'}.}  We retain the notations of section 4.1. Also, without loss of generality we may and will assume that $G$ is nonconstant and $p,q>0$.

Let us consider the open set
$\Di\times\Di_{r_{0}(t)}(c)\subset\Co^{2}$. By
$v=(x,c+y)\in \overline\Di\times\overline{\Di}_{r_{0}(t)}(c)$, $|x|=1$, $|y|=r_{0}(t)$, we
denote a point such that
$$
|G(v)|=\sup_{\Di\times\Di_{r_{0}(t)}(c)}|G|.
$$
We set 
\begin{equation}\label{s}
s:=\frac{r_{0}(t)e^{-\max\{p,q\}/p}}{12\delta(t)}.
\end{equation}
(Observe that $12\delta(t)> 12$, $1<t\leq 9$, see (\ref{e29'}), (\ref{e23}). Hence, $s<r_{0}(t)/12$.)

Consider the point $v'=(x,c+y+c_{y,s})$ with $c_{y,s}$ as in Theorem \ref{te2}. Then $v'$ belongs to the disk $\{x\}\times\overline{\Di}_{r_{1}}(c)$
of radius $r_{1}:=r_{0}(t)+s$. Applying (\ref{e13'}) to
$h:=G|_{\{x\}\times\Di_{er_{0}(t)}(c)}$ from Hadamard's three circle
inequality (see (\ref{e22})) for disks centered at $(x,c)$ of radii $r_{0}(t)$, $r_{1}$
and $er_{0}(t)$ (observe that $er_{1}<2$, see (\ref{set}), so that $\Di_{er_{0}(t)}(c)\subset\Di_{er_{1}}(c)\subset\Di_{3}$) we obtain that
\begin{equation}\label{e34}
\sup_{\{x\}\times\Di_{r_{1}}(c)}|h|<\left(1+\frac{s}{r_{0}(t)}\right)^{q}
\sup_{\{x\}\times\Di_{r_{0}(t)}(c)}|h|.
\end{equation}
Also, by Theorem \ref{markov} and (\ref{kap1}) with $R:=r_{1}$, $t=e$ we get
\begin{equation}\label{e35}
\sup_{\{x\}\times\Di_{r_{1}}(c)}|h'|\leq\frac{9q}{r_{1}}
\sup_{\{x\}\times\Di_{r_{1}}(c)}|h|.
\end{equation}
Using (\ref{e34}) and (\ref{e35}) we can estimate $|h(v')|$ by the
mean-value inequality:
$$
\begin{array}{c}
\displaystyle |h(v)|-|h(v')|\leq |h(v)-h(v')|\leq
\sup_{\{x\}\times\Di_{r_{1}}(c)}|h'|\cdot
||v-v'||\leq\\
\\
\displaystyle \frac{9qs}{r_{0}(t)+s}\left(1+\frac{s}{r_{0}(t)}\right)^{q}
|h(v)|\leq
\frac{9q}{(e^{\max\{p,q\}/p})(12\delta(t))}
\cdot\left(1+\frac{e^{-\max\{p,q\}/p}}{12\delta(t)}\right)^{q}|h(v)|\\
\\
\displaystyle \leq\frac{9e^{-1}p}{12\delta(t)}e^{(qe^{-\max\{p,q\}/p})/(12\delta(t))}
|h(v)|\leq\frac{p}{3\delta(t)}e^{(e^{-1}p)/(12\delta(t))}
|h(v)|\leq
\\
\\
\displaystyle \frac{p}{3\delta(t)}e^{1/(18e)}|h(v)|< \displaystyle
\frac{2p}{5\delta(t)}|h(v)|<|h(v)|.
\end{array}
$$
(We used here the following inequalities: $\delta(t)\geq N_{F}(t)\geq\frac{1}{\ln\left(\frac{1+t}{2\sqrt{t}}\right)}p\geq\frac{2\sqrt{t}}{(\sqrt{t}-1)^{2}}p\geq\frac{3}{2}p$,
$1<t\leq 9$, see (\ref{e29'}),
$\max\{p,q\}/p\geq 1$ and $xe^{-x}\leq e^{-1}$ for $x\geq 1$.) Hence
\begin{equation}\label{e36}
|h(v')|>\left(1-\frac{2p}{5\delta(t)}\right)|h(v)|.
\end{equation}

Next, let us consider the line $l:=\{(z,w)\in\Co^{2}\ :\
w-c-y-c_{y,s}=0\}$. Then $l$ passes through $v'$ and according to
Theorem \ref{te2} intersects the graph $\Gamma=\{(z,F(z))\ :\
z\in\Di\}$ in a set $Y$ containing at least $N_{F}(t)$ points with
pairwise distances greater than $\frac{s(t-1)}{\sqrt{\lambda(t)}}$. For $R>0$ we
set
$$
l_{R}:=\{(z,w)\in l\ :\ |z|<R\}.
$$
We naturally identify $l_{R}$ with disk $\Di_{R}$. Now let us apply
Theorem \ref{cart2} to the univariate function $g_{l}:=G|_{l}$ with $r=1$,
$H:=\frac{s(t-1)}{4t\sqrt{\lambda(t)}}$. According to this theorem and
inequality (\ref{ca6}) in the disk $l_{1}$ outside the union of
open disks $\{D_{j}\}_{1\leq j\leq k}\subset l$, $k<\frac{1}{\ln\left(\frac{1+t}{2\sqrt{t}}\right)}\cdot p$,
with the sum of radii $\sum_{k=1}^{n} r_{k}\leq 2Ht\leq
\frac{s(t-1)}{2\sqrt{\lambda(t)}}$,
\begin{equation}\label{e37}
|g_{l}(z)|\geq M_{g_{l}}(1)\left(
\frac{M_{g_{l}}(1)}{M_{g_{l}}(t)}\right)^{c(H)}\geq M_{g_{l}}(1)e^{-c(H)p}
\end{equation}
where 
\begin{equation}\label{e37'}
c(H):=\frac{(\sqrt{t}+1)^{4}+18(\sqrt{t}+1)^{2}\ln\left(\frac{e}{H}\right)}{2(t-1)^{2}}.
\end{equation}
\begin{R}\label{kbound}
{\rm We single out that the inequality for $k$ is strict.
Indeed, if $g_{l}\equiv const$, then $k=0$ and the conclusion is obvious because $p>0$. If $g_{l}\not\equiv const$, then according to (\ref{ca6})
$$
k\leq\frac{\ln\left(\frac{M_{g_{l}}(tr)}{M_{g_{l}}(\sqrt{t}r)}\right)}{
\ln\left(\frac{1+t}{2\sqrt{t}}\right)}<
\frac{\ln\left(\frac{M_{g_{l}}(tr)}{M_{g_{l}}(r)}\right)}{
\ln\left(\frac{1+t}{2\sqrt{t}}\right)}\leq\frac{p}{\ln\left(\frac{1+t}{2\sqrt{t}}\right)}.
$$
}
\end{R}

Now, by the definition of $H$ we have, see (\ref{set}), (\ref{e29'}),
$$
\begin{array}{c}
\displaystyle
\frac{1}{H}:=\frac{4t\sqrt{\lambda(t)}}{s(t-1)}=\frac{48t\delta(t)\sqrt{\lambda(t)}}{t-1}\left(\frac{4(\sqrt{t}+1)}{M(t-1)}\right)^{\gamma(t)+1}e^{\max\{p,q\}/p}=\\
\\
\displaystyle
\frac{648t(\sqrt{t}+1)^{3}}{(t-1)^{4}}\left(\frac{4(\sqrt{t}+1)}{M(t-1)}\right)^{\gamma(t)+1}\ln\left(\frac{4(\sqrt{t}+1)}{M(t-1)}\right)\sqrt{\ln\left(\frac{2(\sqrt{t}+1)}{M(t-1)^{2}}\right)}\cdot e^{\max\{p,q\}/p}
\\
\\
\displaystyle
<\frac{648t(\sqrt{t}+1)^{3}}{(t-1)^{4}}\left(\frac{32(\sqrt{t}+1)}{M(t-1)^{2}}\right)^{\gamma(t)+2}e^{\max\{p,q\}/p}<\left(\frac{32(\sqrt{t}+1)}{M(t-1)^{2}}\right)^{\gamma(t)+5}e^{\max\{p,q\}/p}.
\end{array}
$$
(We used here that $\frac{4(\sqrt{t}+1)}{M(t-1)}> 10$, $1<t\leq 9$, see (\ref{e23}), $\ln x<\sqrt{x}$ for $x\geq 10$ and
$\ln x<x$ for all $x>1$.)

Hence, see (\ref{set}),
\begin{equation}\label{eq417}
\begin{array}{c}
\displaystyle
e^{c(H)}=e^{\tilde a_{1}(t)}\left(\frac{32(\sqrt{t}+1)}{M(t-1)^{2}}\right)^{\tilde a_{2}(t)\tilde a_{3}(t)}e^{\frac{\tilde a_{3}(t)(\max\{p,q\}+p)}{p}},\ \ \ {\rm where}\\
\\
\displaystyle
\tilde a_{1}(t):=\frac{(\sqrt{t}+1)^{4}}{2(t-1)^{2}},\ \ \
\tilde a_{2}(t):=\frac{(\sqrt{t}+1)^{4}+18(\sqrt{t}+1)^{2}\ln\left(\frac{4et^{3/2}}{\sqrt{t}-1}\right)}{2(t-1)^{2}}+5,\\
\\
\displaystyle
\tilde a_{3}(t):=\frac{9(\sqrt{t}+1)^{2}}{(t-1)^{2}}.
\end{array}
\end{equation}

Since by our assumption $p\leq\ln\left(\frac{1+t}{2\sqrt{t}}\right)\cdot N_{F}(t)$, there exists a point
$a\in Y$ such that $a\not\in\cup_{j}D_{j}$. Actually, by our choice of
$H$ and the definition of $Y$ we obtain that every $D_{j}$ can contain at most one point of
$Y$. But $\#\{Y\}\geq N_{F}(t)$ and the number of the disks $k<\frac{1}{\ln\left(\frac{1+t}{2\sqrt{t}}\right)}\cdot p\leq N_{F}(t)$. This gives the required result.

From (\ref{e37}) we obtain
\begin{equation}\label{e39}
\sup_{z\in\Di}|G_{F}(z)|\geq |g_{l}(a)|\geq M_{g_{l}}(1)e^{-c(H)p}.
\end{equation}
Using that $M_{g_{l}}(1)\geq |h(v')|$, $\tilde a_{2}(t)>10$, $\tilde a_{3}(t)>1$, $\frac{\tilde a_{1}(t)}{\tilde a_{3}(t)}< 1$, and
$\frac{2p}{5\delta(t)}\leq\frac{4}{15}$, $\delta(t)\geq 1$ for $1<t\leq 9$, we obtain directly from
(\ref{e39}), (\ref{eq417}) and (\ref{e36}) by the choice of $v$:
\begin{equation}\label{e310}
\begin{array}{c}
\displaystyle
\sup_{\Di\times\Di_{r_{0}(t)}(c)}|G|\leq
\frac{e^{c(H)p}}{1-(2p/(5\delta(t))}M_{G_{F}}(1)\leq
e^{(c(H)+6/11)p}M_{G_{F}}(1)<\\
\\
\displaystyle
e^{\tilde a_{3}(t)q}\left(\frac{50(\sqrt{t}+1)}{M(t-1)^{2}}\right)^{\tilde a_{2}(t)\tilde a_{3}(t)p} M_{G_{F}}(1).
\end{array}
\end{equation}
(We also used the inequality $-\ln(1-s)\leq\frac{15}{11}s$,
$0\leq s\leq\frac{4}{15}$, applied to $s:=\frac{2p}{5\delta(t)}$.)

Finally let us consider the polydisk $D:=\Di\times\Di_{2}(c)$. Since
$|c|<1$, $\Di\times\Di\subset D\subset\Di\times\Di_{3}$. Then from
the fact that the Bernstein index of $G$ over vertical disks in $D$
is $\leq q$ we easily obtain
$$
\sup_{\Di\times\Di_{2\cdot e^{-j}}(c)}|G|\leq
e^{q}\sup_{\Di\times\Di_{2\cdot e^{-j-1}}(c)}|G|,\ \ \ j=0,1,\dots,
\lfloor \ln(2/r_{0}(t))\rfloor .
$$
From here, (\ref{set}) and (\ref{e310}) we deduce that
\begin{equation}\label{e311}
\sup_{\Di\times\Di}|G|\leq e^{2q}\left(\frac{4(\sqrt{t}+1)}{M(t-1)}\right)^{q(\gamma(t)+1)}
\sup_{\Di\times\Di_{r_{0}(t)}(c)}|G|\leq
\displaystyle 
\left(\frac{c_{1}(t)}{M}\right)^{c_{2}(t)(p+q)}M_{G_{F}}(1)
\end{equation}
where
\begin{equation}\label{e3112}
c_{1}(t):=\frac{50(\sqrt{t}+1)}{(t-1)^{2}},\ \ \
c_{2}(t):=\frac{18(\sqrt{t}+1)^{2}+81\ln\left(\frac{108e}{\sqrt{t}-1}\right)}{(\sqrt{t}-1)^{4}}.
\end{equation}

(We used that $\tilde a_{3}(t)+2<\gamma(t)+1$ and $\max\{\gamma(t)+1,\tilde a_{2}(t)\tilde a_{3}(t)\}< c_{2}(t)$ for $1<t\leq 9$.)

The proof of Theorem \ref{te1'} is complete.\ \ \ \ \ $\Box$
\sect{\hspace*{-1em}. Proof of Theorem \ref{te1}: Case
$n\geq 2$}
{\bf 6.1.} In the proof we use the following estimate. 
\begin{Lm}\label{le41}
Let $h$ be a nonconstant holomorphic function in the disk ${\Di}_{R/\sqrt{t}}$, $R>0$, $t>1$. Then for each $s\in [t,\infty)$,
$$
v_{h}(R/\sqrt{t})\geq
\max\left\{\frac{1}{k(t,s)}\ln\left(\frac{M_{h}(R/t)}{M_{h}(R/s)}\right),\ 1\right\}
$$
where
$$
k(t,s):=\ln\left(\frac{8e^{\pi^{2}}s\sqrt{t}}{(\sqrt{t}-1)^{2}}\right).
$$
\end{Lm}
{\bf Proof.} We make use of the following result proved in [JO]:

{\em Let $h(z)=\sum_{k=0}^{\infty} a_{k}z^{k} $ be a $p$-valent
holomorphic function in the disk $\Di_{s_{1}}$, $\mu_{p}(s_{1})=\max_{0\leq
k\leq p}|a_{k}|s_{1}^{k}$ and $0<s_{2}<s_{1}$. Then
$$
M_{h}(s_{2})\leq A(p)\mu_{p}(s_{1})(1-s_{2}/s_{1})^{-2p}
$$
where $A(p)=(p+2)2^{3p-1}e^{p\pi^{2}+12}$.}

Applying this result to a function $h$ with $s_{1}=R/\sqrt{t}$, $s_{2}=R/t$ and
$p=v_{h}(R/\sqrt{t})$ from the Cauchy inequality for coefficients of the
Taylor series of $h$ we obtain 
\begin{equation}\label{e50}
\begin{array}{c}
\displaystyle
M_{h}(R/t)\leq A(p)\mu_{p}(R/\sqrt{t})\left(\frac{\sqrt{t}}{\sqrt{t}-1}\right)^{2p}\leq A(p)\mu_{p}(R/s)\left(\frac{s}{\sqrt{t}}\right)^{p}\left(\frac{\sqrt{t}}{\sqrt{t}-1}\right)^{2p}\\
\\
\displaystyle
\leq A(p)\left(\frac{s\sqrt{t}}{(\sqrt{t}-1)^{2}}\right)^{p}
M_{h}(R/s),\ \ \ s\in[t,\infty).
\end{array}
\end{equation}

Apply now (\ref{e50}) to the functions $h^{k}$, $k\in\N$. Since $$
v_{h^{k}}(R/\sqrt{t})\leq kp,\ \ \ M_{h^{k}}(R/t)=(M_{h}(R/t))^{k},\ \ \ M_{h^{k}}(R/s)=(M_{h}(R/s))^{k},
$$
inequality (\ref{e50}) in this case implies
$$
M_{h}(R/t)\leq \left(\lim_{k\to\infty}(A(kp))^{1/k}\right)\left(\frac{s\sqrt{t}}{(\sqrt{t}-1)^{2}}\right)^{p}
M_{h}(R/s)\leq \left(\frac{8e^{\pi^{2}}s\sqrt{t}}{(\sqrt{t}-1)^{2}}\right)^{p}M_{h}(R/s).
$$

This gives the required inequality.\ \ \ \ \ $\Box$\\
\\
{\bf 6.2.} Let $f$ and $g$ satisfy conditions of Theorem \ref{te1} for $n\geq 2$. As in section 5.2 we will assume without loss of generality that $r=1$, $M_{1}=M$ and $M_{2}=1$.
 
By $l_{v}\in {\cal L}_{n}$  we denote the complex line $\{vz\in\Co^{n}\ :\ z\in\Co,\ v\in\Co^{n},\ ||v||=1\}$. We set 
\begin{equation}\label{e41}
f_{v}(z):=f(vz),\ \ \  g_{v}(z,w):=g(vz,w),\ \ \ z\in\Di_{t},\
w\in\Di_{3}.
\end{equation}
(In notations of section 1.3  $f_{v}:=f_{l_{v}}$ and $g_{v}:=g_{l_{v}}$.) According to assumptions of Theorem \ref{te1}, every $g_{v}\in {\cal F}_{p,q}(1;t;1)$ with $p\leq \ln\left(\frac{1+t}{2\sqrt{t}}\right) N_{f}(1,t)$.

Let $v_{*}\in\Co^{n}$, $||v_{*}||=1$,  be such that
\begin{equation}\label{e42}
M_{f}(1/t)=M_{f_{v_{*}}}(1/t).
\end{equation}
\begin{Lm}\label{le42}
Suppose that
\begin{equation}\label{e43}
||v-v_{*}||\leq\gamma:=\frac{(\sqrt{t}-1)\ln t}{9e\sqrt{t}\max\{t,\ln(1/M)\}}.
\end{equation}
Then
$$
M_{f_{v}}(1/s)\leq M_{f}(1/s),\ \ \ s\in [t,\infty),\ \ \
M_{f_{v}}(1/t)\geq \frac{1}{\sqrt{t}}M_{f}(1/t)\geq\frac{1}{\sqrt{t}}M.
$$
\end{Lm}
Observe that from condition (\ref{e12}) similarly to inequality (\ref{e23}) we deduce
\begin{equation}\label{e44}
\frac{1}{M}\geq t^{2}.
\end{equation}
{\bf Proof.} The first inequality is obvious. So, let us prove the second one.

Due to inequalities (\ref{e32}), (\ref{case1}) we have (recall that $1<t\leq 9$)
\begin{equation}\label{e46}
M_{D_{s}f}(1/t)\leq\frac{9e}{\ln t}\max\{t,\ln(1/M)\}M_{f}(1/t)\ \ \ {\rm for\ all}\ \ \ s\in\Co^{n},\ ||s||=1.
\end{equation}

From (\ref{e46}) by the mean-value inequality we get for
$|z|\leq 1/t$, 
$$
|f_{v}(z)-f_{v_{*}}(z)|\leq\frac{9e}{\ln t}\max\{t,\ln(1/M)\}M_{f}(1/t)||v-v_{*}||\leq
\frac{\sqrt{t}-1}{\sqrt{t}}M_{f}(1/t).
$$
This and (\ref{e42}) imply the required inequality of the lemma:
$$
M_{f_{v}}(1/t)\geq \frac{1}{\sqrt{t}}M_{f}(1/t)\geq\frac{1}{\sqrt{t}}M.\ \ \ \ \ \Box
$$

Thus for each $v$ satisfying (\ref{e43}) we get from Lemmas \ref{le41}, \ref{le42} (see (\ref{e10}), (\ref{e11})): 
\begin{equation}\label{e47}
R_{f_{v}}(1,t,t^{2}):=\frac{M_{f_{v}}(1/t)}{M_{f_{v}}(1/t^{2})}\geq\frac{1}{\sqrt{t}}R_{f}(1,t,t^{2})\geq\sqrt{t}\ \ \ {\rm and}
\end{equation}
\begin{equation}\label{e48}
\begin{array}{c}
\displaystyle
N_{f_{v}}(1,t):=v_{f_{v}}(1/\sqrt{t})\geq\max\left\{\sup_{s\in [t,\infty)}\left\{\frac{\ln R_{f_{v}}(1,t,s)}{k(t,s)}\right\},V_{f}(1,t)\right\}\geq\\
\\
\displaystyle
\max\left\{\sup_{s\in [t,\infty)}\left\{\frac{\ln(R_{f}(1,t,s)/\sqrt{t})}{k(t,s)}\right\}, V_{f}(1,t)\right\}=:N_{f}(1,t).
\end{array}
\end{equation}

Inequalities (\ref{e47}), (\ref{e48}) show that we can apply to the functions $f_{v}$ and $g_{v}\in {\cal F}_{p,q}(1;t;1)$, $p\leq\ln\left(\frac{1+t}{2\sqrt{t}}\right)\cdot N_{f}(1,t)$, the inequality of Theorem \ref{te1} for $n=1$ (proved in section 4) with
$M_{1}:=\frac{1}{\sqrt{t}}M$, $M_{2}:=1$. Then we obtain
\begin{Proposition}\label{pr41}
For $v$ satisfying (\ref{e43}) and $(g_{f})_{v}:=
g_{v}(z,f_{v}(z))$
the inequality
$$
\sup_{\Di\times\Di}|g_{v}|\leq\left(\frac{c_{1}(t)\sqrt{t}}{M}\right)^{c_{2}(t)(p+q)}
M_{(g_{f})_{v}}(1)
$$
holds with $c_{1}(t), c_{2}(t)$ defined by (\ref{c1c2}).\ \ \ \ \  $\Box$
\end{Proposition}

Let $K\subset\Bo^{n}$ be the convex body determined by the formula
\begin{equation}\label{e49}
K:=\{zv\in\Bo^{n}\ :\ (z,v)\in\Co\times\Co^{n},\ 
|z|<1,\ ||v||=1,\ ||v-v_{*}||\leq\gamma\}.
\end{equation}

Then from Proposition \ref{pr41} we obtain the following statement:

{\em Under the assumptions of Theorem \ref{te1},}
\begin{equation}\label{e410}
\sup_{K\times\Di}|g|\leq\left(\frac{c_{1}(t)\sqrt{t}}{M}\right)^{c_{2}(t)(p+q)}
\sup_{K}|g_{f}|.
\end{equation}

We deduce from here the required inequality of the theorem.
\begin{Lm}\label{le44}
For every boundary point $z$ of $\Bo^{n}$
there is a real straight line $l_{z}\subset\Co^{n}$
passing through $x$ such that $l_{z}$ intersects
$K$ in an interval $I_{z}$ of length 
$$
|I_{z}|>\gamma_{1}:=\frac{9}{10}\gamma.
$$
\end{Lm}
{\bf Proof.} Without loss of generality we may assume that
$z$ does not belong to $K$ (for otherwise, choose $l_{z}$ joining $z$ with $0$ so that $|I_{z}|=1>\gamma$). Then
as the $l_{z}$ we will take the line passing through
$z$ and $v_{*}$. Considering the real two-dimensional plane $P$
containing $0,z$ and $v_{*}$ we reduce the question on the bound of $|I_{z}|$ to the two-dimensional case. Without loss of generality we may assume that $P=\Re^{2}$ and $v_{*}=(1,0)\in\Re^{2}$. In this case
$K_{P}:=K\cap P$ is a convex set defined in polar coordinates $(r,\phi)$ by the inequalities
$$
|2\sin(\phi/2)|\leq\gamma,\ \ \ |r|<1.
$$

Also, we may assume that $z$ belongs to the upper semicircle $S_{+}$ of the unit disk. Now, from (\ref{e43}) and the inequality $1<t\leq 9$ we obtain that $K_{P}\cap S_{+}\subset\{(r,\phi)\ :\ 0\leq\phi<\pi/3\}$. Thus $|I_{z}|$ is $\geq$ the distance from $v_{*}$ to the line $\{(r,\phi)\ :\ \sin(\phi/2)=\gamma/2\}$ that equals $\sin(2\sin^{-1}(\gamma/2)):=\gamma\cdot\sqrt{1-(\gamma/2)^{2}}$. From here and
(\ref{e43}) we obtain
$$
|I_{z}|>\gamma\sqrt{1-\left(\frac{2\ln 3}{9e}\right)^{2}}>\frac{9}{10}\gamma.\ \ \ \ \
\Box
$$

Let $(x,y)$ be the boundary point of $\Bo^{n}\times\Di$ such that 
\begin{equation}\label{e411}
|g(x,y)|=\sup_{\Bo^{n}\times\Di}|g|.
\end{equation}
According to Lemma \ref{le44} there is a straight line $l\subset\Co^{n}\times\{y\}$ passing trough $(x,y)$ and intersecting $K\times\{y\}$ in the interval $I$ of length  
$>\gamma_{1}$. Let $l^{c}$ be the complex line containing $l$.
We set 
$$
D_{1}:=l^{c}\cap (\Bo^{n}\times\{y\})\ \ \ {\rm and} \ \ \
D_{2}:=l^{c}\cap (\Bo_{t}^{n}\times\{y\}).
$$
We can naturally identify $D_{1}$ and $D_{2}$ we the disks
centered at the point $o\in l^{c}$ such that $d:=||o-(0,y)||:=dist(0,l^{c}-(0,y))$ of radii $r_{1}:=\sqrt{1-d^{2}}$ and $r_{2}:=\sqrt{t^{2}-d^{2}}$. Observe also that
$t r_{1}\leq r_{2}$. Thus $tD_{1}\subset D_{2}$ where
$kD_{1}$, $k>0$, denotes the dilation of $D_{1}$ in $k$ times with respect to $o$. 

Further, by (\ref{e411}) we have for
$\tilde g:=g|_{D_{2}}$:
\begin{equation}\label{e412}
M_{\tilde g}(r_{1};o)=M_{g(\cdot,y)}(1)\ \ \ {\rm and}\ \ \
M_{\tilde g}(tr_{1};o)\leq M_{\tilde g}(r_{2};o)\leq M_{g(\cdot, y)}(t).
\end{equation}
Also, the first condition in (\ref{e13}) implies easily that
\begin{equation}\label{e413}
M_{g(\cdot, y)}(t)\leq e^{p}\cdot M_{g(\cdot,y)}(1).
\end{equation}

By the definition the interval $I$ of length $>\gamma_{1}$ is contained in $D_{1}$. Assuming without loss of generality that
$\tilde g$ is nonconstant, we apply to the triple $\tilde g$, $D_{1}$, $tD_{1}$ Theorem \ref{cart2} with 
$H:=\gamma_{1}/4t$. According to this theorem in the disk $D_{1}$ outside the union of a finite number of disks with the sum of radii $<(\gamma_{1}r_{1})/2$ we have
\begin{equation}\label{e414}
\frac{|\tilde g(z)|}{M_{\tilde g}(r_{1};o)}>\left(\frac{M_{\tilde g}(r_{1};o)}{M_{\tilde g}(tr_{1};o)}\right)^{c(H)}
\end{equation}
where 
$$
c(H):=\frac{(\sqrt{t}+1)^{4}+18(\sqrt{t}+1)^{2}\ln\left(\frac{4et}{\gamma_{1}}\right)}{2(t-1)^{2}}.
$$
Since $r_{1}\leq 1$ and $|I|>\gamma_{1}$, the union of such disks cannot cover $I$. In particular, there is a point
$a\in I$ at which inequality (\ref{e414}) holds. Now, from
this inequality and (\ref{e411}), (\ref{e412}), (\ref{e413}),
(\ref{e410}) we obtain
\begin{equation}\label{e415}
\sup_{\Bo^{n}\times\Di}|g|\leq e^{c(H)p}\sup_{K\times\Di}|g|\leq e^{c(H)p}\left(\frac{c_{1}(t)\sqrt{t}}{M}\right)^{c_{2}(t)(p+q)}M_{g_{f}}(1).
\end{equation}
Observe that
\begin{equation}\label{eq516}
\begin{array}{c}
\displaystyle
e^{c(H)}=e^{\tilde a_{1}(t)}\cdot\left(\frac{40e^{2}t^{3/2}\max\{t,\ln(1/M)\}}{(\sqrt{t}-1)\ln t}\right)^{\tilde a_{3}(t)}\ \ \
{\rm where}\\
\\
\displaystyle
\tilde a_{1}(t)=\frac{(\sqrt{t}+1)^{4}}{2(t-1)^{2}},\ \ \
\tilde a_{3}(t)=\frac{9(\sqrt{t}+1)^{2}}{(t-1)^{2}}.
\end{array}
\end{equation}
Also, recall that
$$
c_{1}(t):=\frac{50(\sqrt{t}+1)}{(t-1)^{2}},\ \ \
c_{2}(t):=\frac{18(\sqrt{t}+1)^{2}+81\ln\left(\frac{108e}{\sqrt{t}-1}\right)}{(\sqrt{t}-1)^{4}}.
$$
We consider two cases:

(1) If $\max\{t,\ln(1/M)\}=t$, then (because $1/M\geq t^{2}$)
$$
\begin{array}{c}
\displaystyle
\frac{40e^{2}t^{3/2}t}{(\sqrt{t}-1)\ln t}\leq\frac{40e^{2}(\sqrt{t}+1)\sqrt{t}}{M(t-1)\ln t}=
\frac{4e^{2}}{5}\cdot\frac{t-1}{t\ln t}\cdot \frac{c_{1}(t)t^{3/2}}{M}\leq\\
\\
\displaystyle
\frac{4e^{2}}{5}\cdot 1\cdot \frac{c_{1}(t)t^{3/2}}{M}<\frac{6c_{1}(t)t^{3/2}}{M}.
\end{array}
$$
(We used that the function $t\mapsto\frac{t-1}{t\ln t}$, $t> 1$, is decreasing.)

Also, since $c_{1}(t)>e$, for $1<t\leq 9$ we have
$$
e^{\tilde a_{1}(t)+1}<(c_{1}(t)t^{3/2})^{\tilde a_{1}(t)+1}<(c_{1}(t)t^{3/2})^{\frac{81\ln\left(\frac{108e}{\sqrt{t}-1}\right)}{(\sqrt{t}-1)^{4}}}.
$$
Combining together these inequalities we get in this case
\begin{equation}\label{eq517}
e^{c(H)+1}\leq\left(\frac{6c_{1}(t)t^{3/2}}{M}\right)^{c_{2}(t)}.
\end{equation}

(2) Assume now that $\max\{t,\ln(1/M)\}=\ln(1/M)$. Then
we obtain as before
$$
\begin{array}{c}
\displaystyle
\frac{40e^{2}t^{3/2}\max\{t,\ln(1/M)\}}{(\sqrt{t}-1)\ln t}=
\frac{4e^{2}}{5}\cdot M\ln(1/M)\cdot\frac{t-1}{\ln t}\cdot\frac{c_{1}(t)t^{3/2}}{M}\leq\\
\\
\displaystyle
\frac{4e^{2}}{5}\cdot t^{2}e^{-t}\cdot\frac{t-1}{t\ln t}\cdot
\frac{c_{1}(t)t^{3/2}}{M}\leq \frac{16}{5}\cdot\frac{c_{1}(t)t^{3/2}}{M}.
\end{array}
$$
(We used that $M\ln(1/M)\leq te^{-t}$, because $M\leq e^{-t}<e^{-1}$, and that $t^{2}e^{-t}\leq 4e^{-2}$.)

Thus in this case
\begin{equation}\label{eq518}
e^{c(H)+1}\leq\left(\frac{4c_{1}(t)t^{3/2}}{M}\right)^{c_{2}(t)}.
\end{equation}

From (\ref{e415}) and (\ref{eq517}), (\ref{eq518}) we conclude that
\begin{equation}\label{e416}
\sup_{\Bo^{n}\times\Di}|g|\leq
e^{-p}\left(\frac{a_{1}(t)}{M}\right)^{a_{2}(t)(p+q)}M_{g_{f}}(1),
\end{equation}
\begin{equation}\label{e61}
M_{g_{f}}(t)\leq\sup_{\Bo_{t}^{n}\times\Di}|g|
\leq e^{p}\sup_{\Bo^{n}\times\Di}|g|\leq
\left(\frac{a_{1}(t)}{M}\right)^{a_{2}(t)(p+q)}M_{g_{f}}(1)
\end{equation}
where
$$
a_{1}(t):=6c_{1}(t)t^{3/2},\ \ \ a_{2}(t):=2c_{2}(t).
$$

This completes the proof of Theorem \ref{te1}.\ \ \ \ \ $\Box$
\sect{\hspace*{-1em}. Proof of Theorem \ref{te12}}
Choosing a suitable permutation of coordinates on $\Co^{k}$ without loss of generality we may assume that $i_{j}=j$, $1\leq j\leq k$. For a fixed $(w_{2},\dots, w_{k})\in\Di_{3M_{22}}\times\cdots\times\Di_{3M_{k2}}$ from the conditions of the theorem for $p_{0}$ by Theorem \ref{te1} we have:
$$
\begin{array}{c}
\displaystyle
\sup_{(z,w_{1})\in\Di_{r}\times\Di_{M_{12}}}|g(z,w_{1},w_{2},\dots, w_{k})|\leq e^{p_{1}}\sup_{z\in\Di_{r}}|g(z,f_{1}(z),w_{2},\dots,w_{k})|\\
\\
\displaystyle
\sup_{z\in\Di_{tr}}|g(z,f_{1}(z),w_{2},\dots,w_{k})|\leq e^{p_{1}}\sup_{z\in\Di_{r}}|g(z,f_{1}(z),w_{2},\dots,w_{k})|.
\end{array}
$$
Further, for $(w_{3},\dots, w_{k})\in\Di_{3M_{32}}\times\cdots\times\Di_{3M_{k2}}$ from the conditions of the theorem for $p_{1}$ we have by Theorem \ref{te1} and the previous inequality:
$$
\begin{array}{c}
\displaystyle
\sup_{(z,w_{2})\in\Di_{r}\times\Di_{M_{22}}}|g(z,f_{1}(z),w_{2},\dots, w_{k})|\leq e^{p_{2}}\sup_{z\in\Di_{r}}|g(z,f_{1}(z),f_{2}(z),w_{3}\dots,w_{k})|\\
\\
\displaystyle
\sup_{z\in\Di_{tr}}|g(z,f_{1}(z),f_{2}(z),w_{3}\dots,w_{k})|\leq e^{p_{2}}\sup_{z\in\Di_{r}}|g(z,f_{1}(z),f_{2}(z),w_{3}\dots,
w_{k})|.
\end{array}
$$
Continuing this process we finally obtain 
$$
\begin{array}{c}
\displaystyle
\sup_{(z,w_{k})\in\Di_{r}\times\Di_{M_{k2}}}|g(z,f_{1}(z ),\dots,f_{k-1}(z), w_{k})|\leq e^{p_{k}}\sup_{z\in\Di_{r}}|g(z,f_{1}(z),\dots,f_{k}(z))|\\
\\
\displaystyle
\sup_{z\in\Di_{tr}}|g(z,f_{1}(z),\dots,f_{k}(z))|\leq e^{p_{k}}\sup_{z\in\Di_{r}}|g(z,f_{1}(z),\dots,
f_{k}(z))|.
\end{array}
$$
Combining together all previous inequalities we get the required:
$$
\max_{\Di_{r}\times\Di_{M_{12}}\times\cdots\times\Di_{M_{k2}}}|g|\leq e^{p_{1}+\cdots +p_{k}}M_{g_{\Phi}}(r)\ \ \ {\rm and}\ \
\
M_{g_{\Phi}}(tr)\leq e^{p_{k}}M_{g_{\Phi}}(r).\ \ \ \ \ \ \Box
$$
\sect{\hspace*{-1em}. Proof of Theorem \ref{te13}}
{\bf 8.1.} We first prove the theorem for nonpolynomial entire functions on $\Co^{n}$ of order $\rho<\infty$. We use the following result established in [L, Theorem I.16].

Suppose that $\theta(x)$, $x>x_{*}>0$, is a positive function with
$$
\rho=\limsup_{x\to\infty}\frac{\ln\theta(x)}{\ln x}<\infty.
$$
Then $\theta$ has a {\em proximate order} $\rho(x)$ with the following properties:
\begin{itemize}
\item[(i)]
$\lim_{x\to\infty}\rho(x)=\rho$;
\item[(ii)]
$\theta(x)\leq x^{\rho(x)}$ and $\theta(x_{j})=x_{j}^{\rho(x_{j})}$ for some sequence $x_{j}\to\infty$;
\item[(iii)]
the function $\psi(x)=x^{\rho(x)-\rho}$ is {\em slowly increasing}, i.e.,
$$
\lim_{x\to\infty}\frac{\psi(kx)}{\psi(x)}=1
$$
uniformly on each interval $0<a\leq k\leq b<\infty$;
\item[(iv)]
$\rho(x)$ is a smooth function satisfying
$$
\lim_{x\to\infty}x\rho'(x)\ln x=0.
$$
\end{itemize}

In fact, condition (iii) follows from condition (iv).
Also, if $x^{\rho(x)-\rho}$ is a slowly increasing function, then for every $\epsilon>0$ and every $0< a<b<\infty$ there is $x_{0}$ such that
\begin{equation}\label{e71}
(1-\epsilon)k^{\rho}x^{\rho(x)}<(kx)^{\rho(kx)}<(1+\epsilon)k^{\rho}x^{\rho(x)}
\end{equation}
for $a\leq k\leq b$ and $x\geq x_{0}$.
 
Now, let $f$ be a nonpolynomial entire function on $\Co^{n}$ of order $\rho<\infty$.  As before, we define
\begin{equation}\label{e72}
\phi_{f}(t):=m_{f}(e^{t}),\ \ \ t\in\Re,
\end{equation}
where $m_{f}(r):=\ln M_{f}(r)$, $r>0$.
Then $\phi_{f}$ is a convex increasing function. In particular,
$\phi_{f}'$ is a positive nondecreasing function on $\Re$
(here $\phi_{f}'$ is defined before the formulation of Theorem \ref{te14}). Let $\rho_{f}$ be the proximate order of $m_{f}$.
We define
$$
\tilde\rho_{f}(t):=\rho_{f}(e^{t}).
$$
\begin{Lm}\label{l71}
Set
$$
\alpha_{\rho}:=\min\{1,\ln(1+1/\rho)\}.
$$
Then
$$
\liminf_{t\to\infty}\frac{\phi_{f}'(t+\alpha_{\rho})+\rho e^{t\tilde\rho_{f}(t)}}{\phi_{f}'(t-2\alpha_{\rho})}\leq e^{3}+e^{2}.
$$
\end{Lm}
{\bf Proof.} Assume, on the contrary, that the statement of the lemma is wrong. Then there are positive numbers $t_{0}$, $a$ such that for any $t\geq t_{0}$
\begin{equation}\label{el81}
\phi_{f}'(t+\alpha_{\rho})+\rho e^{t\tilde\rho_{f}(t)}>(e^{3}+e^{2}+a)\phi_{f}'(t-2\alpha_{\rho}).
\end{equation}
We will assume without loss of generality that $\rho>0$. (The arguments in the case $\rho=0$ are similar.) Then conditions (i) and (iv) in the definition of $\rho_{f}$ yield
$$
\lim_{t\to\infty}\frac{\rho}{(t\tilde\rho_{f}(t))'}=1.
$$
From here and (\ref{el81}) for sufficiently large $t_{0}$ we obtain
$$
\phi_{f}'(t+\alpha_{\rho})+b (e^{t\tilde\rho_{f}(t)})'>(e^{3}+e^{2}+a)\phi_{f}'(t-2\alpha_{\rho})
$$
where $b:=1+ae^{-3}$.

Integrating this inequality from $t_{0}$ to $t$ we get
\begin{equation}\label{e73}
\begin{array}{c}
\displaystyle
\phi_{f}(t+\alpha_{\rho})-\phi_{f}(t_{0}+\alpha_{\rho})+be^{t\tilde\rho_{f}(t)}-be^{t_{0}\tilde\rho_{f}(t_{0})}>\\
\\
\displaystyle
(e^{3}+e^{2}+a) (\phi_{f}(t-2\alpha_{\rho})-\phi_{f}(t_{0}-2\alpha_{\rho})).
\end{array}
\end{equation}

Next, by $t_{j}:=\ln x_{j}$, $j\in\N$, we denote the sequence from condition (ii) of the definition of the proximate order $\rho_{f}$. Now, for every $\epsilon>0$ there is $t_{\epsilon}>0$ such that for all $t\geq t_{\epsilon}$ we have from (\ref{e71}) 
\begin{equation}\label{e74}
\begin{array}{c}
\displaystyle
(1-\epsilon)e^{3\rho\alpha_{\rho}}e^{(t-2\alpha_{\rho})\tilde\rho_{f}(t-2\alpha_{\rho})}<
e^{(t+\alpha_{\rho})\tilde\rho_{f}(t+\alpha_{\rho})}<(1+\epsilon)e^{3\rho\alpha_{\rho}}e^{(t-2\alpha_{\rho})\tilde\rho_{f}(t-2\alpha_{\rho})}\\
\\
\displaystyle
(1-\epsilon)e^{2\rho\alpha_{\rho}}e^{(t-2\alpha_{\rho})\tilde\rho_{f}(t-2\alpha_{\rho})}<
e^{t\tilde\rho_{f}(t)}<(1+\epsilon)e^{2\rho\alpha_{\rho}}e^{(t-2\alpha_{\rho})\tilde\rho_{f}(t-2\alpha_{\rho})}.
\end{array}
\end{equation}
Moreover, according to condition (ii) for $\rho_{f}$,
$$
\phi_{f}(t)\leq e^{t\tilde\rho_{f}(t)},\ \ \
\phi_{f}(t_{j})=e^{t_{j}\tilde\rho_{f}(t_{j})}.
$$
From here and inequalities (\ref{e74}) for all $j\geq j_{\epsilon}$ we obtain
\begin{equation}\label{e75}
\phi_{f}(t_{j}+3\alpha_{\rho})\leq (1+\epsilon)e^{3\rho\alpha_{\rho}}\phi_{f}(t_{j}).
\end{equation}
Using (\ref{e74}), (\ref{e75}) and condition (iii) for $\rho_{f}$ we derive from (\ref{e73}) with $t:=t_{j}+2\alpha_{\rho}$
\begin{equation}\label{e76}
\begin{array}{c}
\displaystyle
e^{3}+e^{2}+a\leq\\
\\
\displaystyle
\liminf_{j\to\infty}\frac{\phi_{f}(t_{j}+3\alpha_{\rho})-
\phi_{f}(t_{0}+\alpha_{\rho})+be^{(t_{j}+2\alpha_{\rho})\tilde\rho_{f}(t_{j}+2\alpha_{\rho})}-be^{t_{0}\tilde\rho_{f}( t_{0})}}{\phi_{f}(t_{j})-\phi_{f}(t_{0}-2\alpha_{\rho})}=\\
\\
\displaystyle
\liminf_{j\to\infty}\frac{\phi_{f}(t_{j}+3\alpha_{\rho})+be^{(t_{j}+2\alpha_{\rho})\tilde\rho_{f}(t_{j}+2\alpha_{\rho})}}{\phi_{f}(t_{j})}\leq e^{3\rho\alpha_{\rho}}+be^{2\rho\alpha_{\rho}}\leq
e^{3}+e^{2}+ae^{-1}.
\end{array}
\end{equation}
This gives a contradiction.\ \ \ \ \ $\Box$

As a corollary we obtain
\begin{Proposition}\label{pr72}
There is a sequence $\{\tilde r_{j}\}\subset\Re_{+}$ convergent to $\infty$ such that
$$
\begin{array}{c}
\displaystyle
\frac{m_{f}(e^{\alpha_{\rho}}\tilde r_{j})-m_{f}( e^{-2\alpha_{\rho}}\tilde r_{j})+\ln( a_{1}(e^{\alpha_{\rho}}))}{m_{f}(e^{-\alpha_{\rho}}\tilde r_{j})-m_{f}(e^{-2\alpha_{\rho}}\tilde r_{j}/4)-1}\leq 83,\\
\\
\displaystyle
\frac{\rho \tilde r_{j}^{\rho_{f}(\tilde r_{j})}}{m_{f}(e^{-\alpha_{\rho}}\tilde r_{j})-m_{f}(e^{-2\alpha_{\rho}}\tilde r_{j}/4)-1}\leq\frac{28}{\alpha_{\rho}},\ \ \
j\in\N.
\end{array}
$$
Here $a_{1}(t)$ is defined in (\ref{a1a2}).
\end{Proposition}
{\bf Proof.} As the $\tilde r_{j}$ we will take $e^{s_{j}}$ where
$\{s_{j}\}\subset\Re_{+}$ is such that
\begin{equation}\label{e77}
\begin{array}{c}
\displaystyle
\frac{\phi_{f}'(s_{j}+\alpha_{\rho})+\ln (a_{1}(e^{\alpha_{\rho}}))/\alpha_{\rho}}{\phi_{f}'(s_{j}-2\alpha_{\rho})-1/\alpha_{\rho}}\leq e^{3}+e^{2}+1/10,\\
\\
\displaystyle
\frac{\rho e^{s_{j}\tilde\rho_{f}(s_{j}) }}{\phi_{f}'(s_{j}-2\alpha_{\rho})-1/\alpha_{\rho}}\leq e^{3}+e^{2}+1/10,\ \ \ j\in\N.
\end{array}
\end{equation}
Observe that $\phi_{f}'(t)$ tends to $\infty$ as $t\to\infty$.
(For otherwise, $m_{f}(r)\leq A\ln r$ so that $f$ is a polynomial.) Thus such a sequence exists by Lemma \ref{l71}. 

Since the function
$\phi_{f}'$ is nondecreasing from (\ref{e77}) we get
\begin{equation}\label{e78}
\begin{array}{c}
\displaystyle
\frac{m_{f}(e^{\alpha_{\rho}}\tilde r_{j})-m_{f}(e^{-2\alpha_{\rho}}\tilde r_{j})+\ln (a_{1}(e^{\alpha_{\rho}}))}{m_{f}(e^{-\alpha_{\rho}}\tilde r_{j})-m_{f}(e^{-2\alpha_{\rho}}\tilde r_{j})-1}=\\
\\
\displaystyle
\frac{\phi_{f}(s_{j}+\alpha_{\rho})-\phi_{f}(s_{j}-2\alpha_{\rho})+\ln (a_{1}(e^{\alpha_{\rho}}))}{\phi_{f}(s_{j}-\alpha_{\rho})-\phi_{f}(s_{j}-2\alpha_{\rho})-1}\leq
\\
\\
\displaystyle
\frac{(3\alpha_{\rho})\phi_{f}'(s_{j}+\alpha_{\rho})+\ln (a_{1}(e^{\alpha_{\rho}}))}{(\alpha_{\rho})\phi_{f}'(s_{j}-2\alpha_{\rho})-1}\leq
3(e^{3}+e^{2}+1/10)<83.
\end{array}
\end{equation}

Similarly,
\begin{equation}\label{e78'}
\begin{array}{c}
\displaystyle
\frac{\rho \tilde r_{j}^{\rho_{f}(\tilde r_{j})}}{m_{f}(e^{-\alpha_{\rho}}\tilde r_{j})-m_{f}(e^{-2\alpha_{\rho}}\tilde r_{j}/4)-1}\leq\frac{\rho e^{s_{j}\tilde\rho_{f}(s_{j})}}
{(\alpha_{\rho})\phi_{f}'(s_{j}-2\alpha_{\rho})-1}\leq\\
\\
\displaystyle
\frac{e^{3}+e^{2}+1/10}{\alpha_{\rho}}<\frac{28}{\alpha_{\rho}}.
\ \ \ \ \ \Box
\end{array}
\end{equation}

Now, let us prove the theorem in the case $\rho<\infty$.

For the sequence $\{\tilde r_{j}\}$ satisfying Proposition \ref{pr72} we have, see (\ref{e11}) and the definition of $\alpha_{\rho}$,
\begin{equation}\label{e79}
\begin{array}{c}
\displaystyle
N_{f}(\tilde r_{j},e^{\alpha_{\rho}})\geq\frac{m_{f}(e^{-\alpha_{\rho}}\tilde r_{j})-m_{f}(e^{-2\alpha_{\rho}}\tilde r_{j})-\alpha_{\rho}/2}{k(e^{\alpha_{\rho}},e^{2\alpha_{\rho}})}>
\\
\\
\displaystyle
\frac{m_{f}(e^{-\alpha_{\rho}}\tilde r_{j})-m_{f}(e^{-2\alpha_{\rho}}\tilde r_{j})-1}{\ln 8+ \pi^{2}+(5\alpha_{\rho})/2+2\ln(\sqrt{e}+1)+2\ln(\rho+1/(e-1))}>
\\
\\
\displaystyle
\frac{m_{f}(e^{-\alpha_{\rho}}\tilde r_{j})-m_{f}(e^{-2\alpha_{\rho}}\tilde r_{j})-1}{17+2\ln(\rho+1)}.
\end{array}
\end{equation}
Set
\begin{equation}\label{e710}
n_{j}:=\frac{m_{f}(e^{-\alpha_{\rho}}\tilde r_{j})-m_{f}(e^{-2\alpha_{\rho}}\tilde r_{j})-1}{9(\sqrt{e}+1)^{2}(\rho^{2}+1)(17+2\ln(\rho+1))},\ \ \ j\in\N.
\end{equation}
Then using inequality (\ref{ca8}) we obtain
\begin{equation}\label{e711}
n_{j}\leq\ln\left(\frac{1+e^{\alpha_{\rho}}}{2e^{\alpha_{\rho}/2}}\right)N_{f}(\tilde r_{j},e^{\alpha_{\rho}}).
\end{equation}
Also, according to Proposition \ref{pr72} we have
\begin{equation}\label{e712}
m_{f}(e^{\alpha_{\rho}}\tilde r_{j})-m_{f}(e^{-\alpha_{\rho}}\tilde r_{j})+\ln(a_{1}(e^{\alpha_{\rho}}))\leq c(\rho)n_{j},\ \ \
i\in\N,
\end{equation}
where
\begin{equation}\label{e713}
c(\rho):=747(\sqrt{e}+1)^{2}(\rho^{2}+1)(17+2\ln(\rho+1))\leq c
(\rho+1)^{2}(1+\ln(\rho+1))
\end{equation}
for an absolute constant $c>0$.

Observe that $n_{j}\to\infty$ as $j\to\infty$ because $f$ is not a polynomial. In addition, we have the following statement.
\begin{Proposition}\label{pr73}
For all sufficiently large $j$ the inequality
$$
n_{j}^{1/(\rho+\epsilon_{j})}\leq\tilde r_{j}\leq\left(\frac{A}{\rho_{*}}\right)^{1/\rho_{*}}\ \!n_{j}^{1/(\rho-\epsilon_{j}')}
$$
holds for an absolute constant $A>0$ and some sequences $\{\epsilon_{j}\}, \{\epsilon_{j}'\}\subset\Re_{+}$ convergent to $0$. Here
$\rho_{*}:=\min\{1,\rho\}$.
\end{Proposition}
(In the case $\rho=0$ we assume that the right-hand side is $\infty$.)\\
{\bf Proof.} According to formula (\ref{e710}) and the definition of the order $\rho$, for all sufficiently large $j$ the inequality
$$
n_{j}\leq\frac{83m_{f}(\tilde r_{j})}{c(\rho)}\leq\frac{83\tilde r_{j}^{\rho+\tilde\epsilon_{j}}}{c(\rho)}\leq\tilde r_{j}^{\rho+\epsilon_{j}}
$$
holds
for some sequences $\{\tilde \epsilon_{j}\}, \{\epsilon_{j}\}\subset\Re_{+}$ convergent to $0$. This proves the left-hand side inequality of the proposition.

Further, according to Proposition \ref{pr72} and definitions of $c(\rho)$, $\alpha_{\rho}$ and $\rho_{f}$ we have
$$
\tilde r_{j}\leq\left(\frac{28 c(\rho)}{83\rho\alpha_{\rho}}\right)^{1/\rho_{f}(\tilde r_{j})}n_{j}^{1/\rho_{f}(\tilde r_{j})}\leq
\left(\frac{A}{\rho_{*}}\right)^{1/\rho_{*}}n_{j}^{1/(\rho-\epsilon_{j}')}
$$
for an absolute constant $A>0$ and a sequence $\{\epsilon_{j}'\}\subset\Re_{+}$ convergent to $0$.\ \ \ \ \ $\Box$

Suppose now that $g\in {\cal F}_{p,q}(\tilde r_{j};e;M_{f}(e\tilde r_{j}))$
with $p\leq n_{j}$. Then inequality (\ref{e711}) implies that
$g$ satisfies conditions of Theorem \ref{te1} with $r=\tilde r_{j}$, $t=e^{\alpha_{\rho}}$. From this theorem we obtain, see
(\ref{a1a2}),
$$
\begin{array}{c}
\displaystyle
\!\!\sup_{\Bo_{\tilde r_{j}}^{n}\times\Di_{M_{f}(e^{\alpha_{\rho}}\tilde r_{j})}}\!\!\!\!\!\!\!\!\ln |g|\leq
a_{2}(e^{\alpha_{\rho}})(p+q)(\ln a_{1}(e^{\alpha_{\rho}})+m_{f}(e^{\alpha_{\rho}}\tilde r_{j})-m_{f}(e^{-\alpha_{\rho}}\tilde r_{j}))
+\ln M_{g_{f}}(\tilde r_{j}),\\
\\
\displaystyle
\ln M_{g_{f}}(e^{\alpha_{\rho}}\tilde r_{j})
\leq
a_{2}(e^{\alpha_{\rho}})(p+q)(\ln a_{1}(e^{\alpha_{\rho}})+m_{f}(e^{\alpha_{\rho}}\tilde r_{j})-m_{f}(e^{-\alpha_{\rho}}\tilde r_{j}))
+\ln M_{g_{f}}(\tilde r_{j}).
\end{array}
$$
Estimating $a_{2}(e^{\alpha_{\rho}})$ by (\ref{a1a2}) from here, (\ref{e712}) and (\ref{e713})  we deduce that
\begin{equation}\label{e714}
\begin{array}{c}
\displaystyle
\sup_{\Bo_{\tilde r_{j}}^{n}\times\Di_{M_{f}(e^{\alpha_{\rho}}\tilde r_{j})}} |g|\leq e^{C(\rho)n_{j}\max\{p,q\}}M_{g_{f}}(\tilde r_{j}),\\
\\
\displaystyle
M_{g_{f}}(e^{\alpha_{\rho}}\tilde r_{j})\leq e^{C(\rho)n_{j}\max\{p,q\}}M_{g_{f}}(\tilde r_{j}) 
\end{array}
\end{equation}
where $C(\rho)=C(\rho+1)^{6}(1+\ln(\rho+1))^{2}$ for an absolute constant $C>0$.

Finally, as the sequence $\{r_{j}\}$ of the theorem we will take $\{\frac{\tilde r_{j}}{e}\}$ where $\{\tilde r_{j}\}$ satisfies Propositions \ref{pr72} and \ref{pr73}.
Then by the Hadamard three circle inequality, see section 3.1,
we have
$$
\frac{M_{g_{f}}(er_{j})}{M_{g_{f}}(r_{j})}\leq
\frac{M_{g_{f}}(e^{\alpha_{\rho}}\tilde r_{j})}{M_{g_{f}}(e^{\alpha_{\rho}}\tilde r_{j}/e)}\leq
\left(\frac{M_{g_{f}}(e^{\alpha_{\rho}}\tilde r_{j})}{M_{g_{f}}(\tilde r_{j})}\right)^{1/\alpha_{\rho}}
\leq e^{C_{\rho}n_{j}\max\{p,q\}}
$$
where
$$
C_{\rho}:=\frac{C(\rho)}{\alpha_{\rho}}\leq c(\rho+1)^{7}(1+\ln(\rho+1))^{2}
$$
for and absolute constant $c>0$.

From here, the first inequality (\ref{e714}) and the results of sections 3.1, 3.3, see also Remark \ref{re2}, we obtain straightforwardly the required inequalities of the theorem for $\rho<\infty$.  Observe also, that if $\rho>0$, then Proposition \ref{pr73} implies that
$$
n_{j}^{1/(\rho+\epsilon_{j}')}\leq r_{j}\leq\frac{1}{e}\left(\frac{A}{\rho_{*}}\right)^{1/\rho_{*}}n_{j}^{1/(\rho-\tilde\epsilon_{j}')}
$$
for an absolute constant $A>0$ and sequences $\{\epsilon_{j}'\}, \{\tilde\epsilon_{j}'\}\subset\Re_{+}$ convergent to $0$. This gives the inequalities of statements (1) and (2) of the theorem.

Thus the proof of the theorem for $\rho<\infty$ is complete.\\
\\
{\bf 8.2.} Let us prove now the theorem for $\rho=\infty$.
For $t\geq t_{0}$ with a sufficiently large $t_{0}$ we set
\begin{equation}\label{e722}
\psi_{f}(t):=\frac{1}{\ln\phi_{f}(t)}.
\end{equation}
The function $\psi_{f}$ is decreasing and differentiable outside a countable set $S\subset\Re$ and its derivative $\psi_{f}'$ is continuous outside $S$ and has discontinuities of the first kind at the points of $S$. We extend $\psi_{f}'$ to $S$ in the same way as we extended $\phi_{f}'$. Then we have
$$
\psi_{f}'(t)=\frac{-\phi_{f}'(t)}{\phi_{f}(t)(\ln\phi(t))^{2}},\ \ \ \ \ t\geq t_{0}.
$$
\begin{Lm}\label{le73}
$$
\liminf_{t\to\infty}(-t^{2}\psi_{f}'(t))=0.
$$
\end{Lm}
{\bf Proof.} Since the order of $f$ is $\infty$, there are sequences $\{t_{j}\}, \{a_{j}\}\subset\Re_{+}$ convergent to $\infty$ such that 
$$
\ln\phi_{f}(t_{j})\geq a_{j}t_{j},\ \ \ j\in\N.
$$
This implies 
\begin{equation}\label{e723}
\psi_{f}(t_{j})\leq\frac{\epsilon_{j}}{t_{j}},\ \ \ j\in\N,\ \ \ {\rm where}\ \ \
\epsilon_{j}:=\frac{1}{a_{j}}\to\infty\ \ \ {\rm as}\ \ \
j\to\infty.
\end{equation}

Assume, on the contrary, that there is some $a>0$ such that
$$
\liminf_{t\to\infty}(-t^{2}\psi_{f}'(t))> a.
$$
Then there is $t_{*}\geq t_{0}$ such that
\begin{equation}\label{e727}
-\psi_{f}'(t)>\frac{a}{t^{2}}\ \ \ {\rm for\ all}\ \ \ t\in [t_{*},\infty).
\end{equation}
Integrating this inequality from $t\geq t_{*}$ to $\infty$  we get
\begin{equation}\label{e728}
\psi_{f}(t)>\frac{a}{t},\ \ \ t\in [t_{*},\infty). 
\end{equation}
From here and (\ref{e723}) we obtain for all sufficiently large $j$,
$$
\frac{a}{t_{j}}\leq\frac{\epsilon_{j}}{t_{j}}\ \!,
$$
a contradiction.\ \ \ \ \ $\Box$

According to this lemma there is a sequence $\{v_{j}\}\subset\Re_{+}$ convergent to $\infty$ such that
\begin{equation}\label{e729}
\lim_{j\to\infty}v_{j}^{2}\psi_{f}'(v_{j})=0
\end{equation}
and $\psi_{f}'$ is continuous at each $v_{j}$.

For each sufficiently large $j$ by $s_{j}>0$ we denote a number such that
\begin{equation}\label{e730}
\frac{\phi_{f}(v_{j})}{\phi_{f}(v_{j}-2s_{j})}=e.
\end{equation}
Since $\phi_{f}$ is an increasing function and $\phi_{f}'$ is a positive nondecreasing function, from
(\ref{e730}) we obtain
$$
1=\ln\phi_{f}(v_{j})-\ln\phi_{f}(v_{j}-2s_{j})\leq
2s_{j}\frac{\phi_{f}'(v_{j})}{\phi_{f}(v_{j}-2s_{j})}=
2es_{j}\frac{\phi_{f}'(v_{j})}{\phi_{f}(v_{j})}.
$$
Hence,
\begin{equation}\label{e731}
\frac{1}{s_{j}}\leq\frac{2e\phi_{f}'(v_{j})}{\phi_{f}(v_{j})}.
\end{equation}

Set
\begin{equation}\label{e732}
\tilde s_{j}:=\min\{s_{j},1\},\ \ \ j\in\N.
\end{equation}
Using the inequality $\ln(1/s)\leq\frac{2}{\sqrt{s}}-2$,
$s>0$, from (\ref{e731}) we have for all sufficiently large $j$ (for which, in particular, $\ln \phi_{f}(v_{j})>0$)
\begin{equation}\label{e733}
\frac{\ln(1/\tilde s_{j})}{\ln \phi_{f}(v_{j})}\leq\frac{2/\sqrt{s_{j}}}{\ln \phi_{f}(v_{j})}\leq 2\sqrt{2e}
\left(\frac{\phi_{f}'(v_{j})}{\phi_{f}(v_{j})(\ln \phi_{f}(v_{j}))^{2}}\right)^{1/2}=2\sqrt{2e}(-\psi_{f}'(v_{j}))^{1/2}.
\end{equation}

From here and (\ref{e729}) we obtain that there is a sequence
$\{\epsilon_{j}\}\subset\Re_{+}$ convergent to $0$ such that for all sufficiently large $j$
\begin{equation}\label{e734}
\frac{\ln(1/\tilde s_{j})}{\ln \phi_{f}(v_{j})}\leq\frac{\epsilon_{j}}{v_{j}},\ \ \ {\rm  equivalently},\ \ \ 
\frac{1}{\tilde s_{j}}\leq\left( m_{f}(e^{v_{j}})\right)^{\frac{\epsilon_{j}}{v_{j}}}.
\end{equation}

Next, we set
\begin{equation}\label{e735}
t_{j}:=e^{\tilde s_{j}},\ \ \ \tilde r_{j}:=e^{v_{j}-\tilde s_{j}},\ \ \ j\in\N,
\end{equation}
and apply Theorem \ref{te1} to $f$ defined on $\Bo_{t_{j}\tilde r_{j}}^{n}$.

According to this theorem for every function $g\in {\cal F}_{p,q}(\tilde r_{j};t_{j};M_{f}(t_{j}\tilde r_{j}))$ with 
$p\leq\ln\left(\frac{1+t_{j}}{2\sqrt{t_{j}}}\right)N_{f}(\tilde r_{j},t_{j})$ and all sufficiently large $j$ we have
\begin{equation}\label{e736}
\begin{array}{c}
\displaystyle
\sup_{\Bo_{\tilde r_{j}}^{n}\times\Di_{M_{f}(t_{j}\tilde r_{j})}}|g|\leq (a_{1}(t_{j})M_{f}(t_{j}\tilde r_{j}))^{a_{2}(t_{j})(p+q)}M_{g_{f}}(\tilde r_{j});\\
\\
\displaystyle
M_{g_{f}}(t_{j}\tilde r_{j})\leq (a_{1}(t_{j})M_{f}(t_{j}\tilde r_{j}))^{a_{2}(t_{j})(p+q)}M_{g_{f}}(\tilde r_{j}).
\end{array}
\end{equation}
(We used here that $M_{f}(\tilde r_{j}/t_{j})\geq 1$ for all sufficiently large $j$ because $t_{j}\leq e$ and $\{\tilde r_{j}\}$ tends to $\infty$.)

Further, according to formulas (\ref{a1a2}), (\ref{e734})  we have
\begin{equation}\label{e737}
\begin{array}{c}
\displaystyle
a_{1}(t_{j})\leq\frac{300(\sqrt{e}+1)e^{3/2}}{(e^{\tilde s_{j}}-1)^{2}}\leq\frac{3562}{\tilde s_{j}^{2}}
\leq 3562 \left( m_{f}(t_{j}\tilde r_{j})\right)^{\frac{2\epsilon_{j}}{v_{j}}},\\
\\
\displaystyle
a_{2}(t_{j})\leq\frac{36(\sqrt{e}+1)^{2}+162\ln\left(\frac{108e}{e^{\tilde s_{j}/2}-1}\right)}{(e^{\tilde s_{j}/2}-1)^{4}}\leq
\frac{4042+2592\ln\left(\frac{216e}{\tilde s_{j}}\right)}{\tilde s_{j}^{4}}\leq\\
\\
\displaystyle
\frac{20567+2592\ln\left(\frac{1}{\tilde s_{j}}\right)}{\tilde s_{j}^{4}}\leq\frac{20567}{\tilde s_{j}^{5}}\leq
20567\left(m_{f}(t_{j}\tilde r_{j})\right)^{\frac{5\epsilon_{j}}{v_{j}}}.
\end{array}
\end{equation}
(We used that $\ln x\leq x-1$ for $x>0$.)

Thus from (\ref{e736}), (\ref{e737}) we obtain for all sufficiently large $j$
\begin{equation}\label{e738}
\begin{array}{c}
\displaystyle
\ln\left(\frac{\sup_{\Bo_{\tilde r_{j}}^{n}\times\Di_{M_{f}(t_{j}\tilde r_{j})}}|g|}{M_{g_{f}}(\tilde r_{j})}\right)\leq 41134\cdot
(m_{f}(t_{j}\tilde r_{j}))^{1+\frac{5\epsilon_{j}}{v_{j}}}\max\{p,q\},
\\
\\
\displaystyle
\ln\left(\frac{M_{g_{f}}(t_{j}\tilde r_{j})}{M_{g_{f}}(\tilde r_{j})}\right)\leq 41134\cdot ( m_{f}(t_{j}\tilde r_{j}))^{1+\frac{5\epsilon_{j}}{v_{j}}}\max\{p,q\}.
\end{array}
\end{equation}

Let us estimate from below the expression
$\ln\left(\frac{1+t_{j}}{2\sqrt{t_{j}}}\right)N_{f}(\tilde r_{j},t_{j})$.
\begin{Lm}\label{below}
There is a sequence $\{\delta_{j}\}\subset\Re_{+}$ convergent to $0$ such that for all sufficiently large $j$
$$
\ln\left(\frac{1+t_{j}}{2\sqrt{t_{j}}}\right)N_{f}(\tilde r_{j},t_{j})\geq (m_{f}(t_{j}\tilde r_{j}))^{1-\delta_{j}}.
$$
\end{Lm}
{\bf Proof.}
Using (\ref{ca8}) we obtain
\begin{equation}\label{e739}
\frac{1}{\ln\left(\frac{1+t_{j}}{2\sqrt{t_{j}}}\right)}\leq\frac{9(\sqrt{t_{j}}+1)^{2}}{(t_{j}-1)^{2}}\leq\frac{9(\sqrt{e}+1)^{2}}{\tilde s_{j}^{2}}\leq 64( m_{f}(t_{j}\tilde r_{j}))^{\frac{2\epsilon_{j}}{v_{j}}}.
\end{equation}
Also, for all sufficiently large $j$ by (\ref{e11}), (\ref{k}), (\ref{e730}) we have
\begin{equation}\label{e740}
N_{f}(\tilde r_{j},t_{j})\geq\frac{\ln(M_{f}(\tilde r_{j}/t_{j})/(\sqrt{e}M_{f}(1)))}{k(t_{j},\tilde r_{j})}\geq \frac{ m_{f}(t_{j}\tilde r_{j})}{3k(t_{j},\tilde r_{j})}.
\end{equation}
In turn, for all sufficiently large $j$,
\begin{equation}\label{!}
k(t_{j},\tilde r_{j})\leq\ln\left(\frac{8e^{\pi^{2}}\tilde r_{j}\sqrt{e}(\sqrt{e}+1)^{2}}{\tilde s_{j}^{2}}\right)\leq \ln \tilde r_{j}+15+\frac{2\epsilon_{j}}{v_{j}}\ln m_{f}(t_{j}\tilde r_{j}).
\end{equation}

Next, we estimate $\ln \tilde r_{j}$.

Since the function $\phi_{f}'$ is nondecreasing, 
\begin{equation}\label{!!}
\frac{\phi_{f}(v_{j})-\phi_{f}(t_{0})}{(v_{j}-t_{0})\phi_{f}(v_{j})(\ln\phi_{f}(v_{j}))^{2}}\leq
\frac{\phi_{f}'(v_{j})}{\phi_{f}(v_{j})(\ln\phi_{f}(v_{j}))^{2}}:=
-\psi_{f}'(v_{j})\leq\frac{\tilde\epsilon_{j}}{v_{j}^{2}}
\end{equation}
for some $\{\tilde\epsilon_{j}\}\subset\Re_{+}$ convergent to $0$.

Also,
\begin{equation}\label{!!!}
\lim_{j\to\infty}\frac{v_{j}(\phi_{f}(v_{j})-\phi_{f}(t_{0}))}{(v_{j}-t_{0})\phi_{f}(v_{j})}=1.
\end{equation}
From (\ref{!!}) and (\ref{!!!}) for all sufficiently large $j$ and some sequence $\{\epsilon_{j}'\}\subset\Re_{+}$ convergent to $0$ we obtain
\begin{equation}\label{tilde}
\ln\tilde r_{j}:=v_{j}-\tilde s_{j}<v_{j}\leq \epsilon_{j}'(\ln\phi_{f}(v_{j}))^{2}=\epsilon_{j}'(\ln m_{f}(t_{j}\tilde r_{j}))^{2}.
\end{equation}
Hence, for all sufficiently large $j$
$$
3k(t_{j},\tilde r_{j})\leq \epsilon_{j}''(\ln m_{f}(t_{j}\tilde r_{j}))^{2}
$$
for some $\{\epsilon_{j}''\}\subset\Re_{+}$ convergent to $0$.

From here, (\ref{e740}) and (\ref{e739}) we obtain for all sufficiently large $j$
$$
\ln\left(\frac{1+t_{j}}{2\sqrt{t_{j}}}\right)N_{f}(\tilde r_{j},t_{j})\geq
(m_{f}(t_{j}\tilde r_{j}))^{1-\delta_{j}}
$$
for some $\{\delta_{j}\}\subset\Re_{+}$ convergent to $0$. \ \ \ \ \ $\Box$

We set
\begin{equation}\label{nj}
n_{j}:=(m_{f}(t_{j}\tilde r_{j}))^{1-\delta_{j}}
\end{equation}
with $\{\delta_{j}\}$ satisfying Lemma \ref{below}.

Then from (\ref{e738}) we get for every $g\in {\cal F}_{p,q}(\tilde r_{j};t_{j};M_{f}(t_{j}\tilde r_{j}))$ with $p\leq n_{j}$
\begin{equation}\label{e742}
\begin{array}{c}
\displaystyle
\sup_{\Bo_{\tilde r_{j}}^{n}\times\Di_{M_{f}(t_{j}\tilde r_{j})}}|g|\leq
e^{n_{j}^{1+\tilde\epsilon_{j}}\max\{p,q\}}M_{g_{f}}(\tilde r_{j}),
\\
\\
\displaystyle
M_{g_{f}}(t_{j}\tilde r_{j})\leq  e^{n_{j}^{1+\tilde\epsilon_{j}}\max\{p,q\}}M_{g_{f}}(\tilde r_{j}).
\end{array}
\end{equation}
for some $\{\tilde\epsilon_{j}\}$ convergent to $0$. 

Finally, we define
$$
r_{j}:=\frac{t_{j}\tilde r_{j}}{e^{2}}.
$$
Observe that ${\cal F}_{p,q}(er_{j};e;M_{f}(e^{2}r_{j}))\subset
{\cal F}_{p,q}(\tilde r_{j};t_{j};M_{f}(t_{j}\tilde r_{j}))$,
because $t_{j}\leq e$. Then by
the Hadamard three circle inequality, see section 3.1, using the estimate for $1/\ln t_{j}$, see (\ref{e734}), we have
for all sufficiently large $j$ and $g\in {\cal F}_{p,q}(er_{j};e;M_{f}(e^{2}r_{j}))$, $p\leq n_{j}$,
$$
\frac{M_{g_{f}}(er_{j})}{M_{g_{f}}(r_{j})}\leq
\frac{M_{g_{f}}(t_{j}\tilde r_{j})}{M_{g_{f}}(t_{j}\tilde r_{j}/e)}\leq\left(\frac{M_{g_{f}}(t_{j}\tilde r_{j})}{M_{g_{f}}(\tilde r_{j})}\right)^{\frac{1}{\ln t_{j}}}\leq
e^{n_{j}^{1+\epsilon_{j}}\max\{p,q\}}
$$
for some $\{\epsilon_{j}\}\subset\Re_{+}$ convergent to $0$.

From here, the first inequality (\ref{e742}), definition (\ref{nj}) and the results of sections 3.1, 3.3, see also Remark \ref{re2}, we obtain straightforwardly the required inequalities of the theorem for $\rho=\infty$.

The proof of Theorem \ref{te13} is complete.\ \ \ \ \ $\Box$
\sect{\hspace*{-1em}. Proofs of Theorem \ref{te14} and Corollary \ref{c14}}
{\bf 9.1. Proof of Theorem \ref{te14}.}

(1) Assume that an entire function $f$ on $\Co^{n}$ of order 
$\rho<\infty$ satisfies
\begin{equation}\label{e743}
\limsup_{t\to\infty}\frac{m_{f}(e^{\alpha_{\rho}}r)-m_{f}(e^{-\alpha_{\rho}}r)+\rho e^{\rho t}}{m_{f}(e^{-\alpha_{\rho}}r)-m_{f}(e^{-2\alpha_{\rho}}r)}<A<\infty
\end{equation}
where $\alpha_{\rho}:=\min\{1,\ln(1+1/\rho)\}$. 

Then for a sufficiently large $r_{0}$ and all $r\geq r_{0}$ we have
\begin{equation}\label{e746}
\begin{array}{c}
\displaystyle
\frac{m_{f}(e^{\alpha_{\rho}}r)-m_{f}(e^{-\alpha_{\rho}}r)+\ln (a_{1}(e^{\alpha_{\rho}}))}{m_{f}(e^{-\alpha_{\rho}}r)-m_{f}(e^{-2\alpha_{\rho}}r)-1}<A,\\
\\
\displaystyle
\frac{\rho e^{\rho t}}{m_{f}(e^{-\alpha_{\rho}}r)-m_{f}(e^{-2\alpha_{\rho}}r)-1}<A.
\end{array}
\end{equation}

Similarly to the definition of $n_{j}$ in section 8.1, see (\ref{e710}), we determine
\begin{equation}\label{e747}
k(r):=\frac{m_{f}(e^{-\alpha_{\rho}} r)-m_{f}(e^{-2\alpha_{\rho}} r)-1}{9(\sqrt{e}+1)^{2}(\rho^{2}+1)(17+2\ln(\rho+1))},\ \ \ r\geq r_{0}.
\end{equation}
Then
\begin{equation}\label{e747'}
k(r)\leq\ln\left(\frac{1+e^{\alpha_{\rho}}}{2e^{\alpha_{\rho}/2}}\right)N_{f}(r,e^{\alpha_{\rho}}).
\end{equation}
Also, according to (\ref{e746})
\begin{equation}\label{e747''}
m_{f}(e^{\alpha_{\rho}}r)-m_{f}(e^{-\alpha_{\rho}}r)+\ln (a_{1}(e^{\alpha_{\rho}}))\leq A\tilde c(\rho)k(r),\ \ \  r\geq r_{0},
\end{equation}
where
$\tilde c(\rho):=9(\sqrt{e}+1)^{2}(\rho^{2}+1)(17+2\ln(\rho+1))$.

Further, as in the proof of Proposition \ref{pr73}, using the definition of the order $\rho$, by (\ref{e746}) we obtain
for a sufficiently large $r_{0}$ and all $r\geq r_{0}$,
\begin{equation}\label{e746'}
(k(r))^{1/(\rho+\epsilon'(r))}\leq r\leq\left(\frac{ \tilde c}{\rho_{*}}\right)^{1/\rho_{*}}(k(r))^{1/\rho},
\end{equation}
for some $\tilde c$ depending on $A$; here $\rho_{*}:=\min\{1,\rho\}$ and $\epsilon':[r_{0},\infty)\to\Re_{*}$ is a continuous function
decreasing to $0$ as $r\to\infty$.

As in section 8.1 inequalities (\ref{e747'}), (\ref{e747''}), (\ref{e746'}) imply the fulfillment of Theorem \ref{te13} for functions $g\in {\cal F}_{p,q}(er;e;M_{f}(e^{2}r))$ with $p\leq k(r)$ in which $n_{j}$ is substituted for $k(r)$, $r\geq r_{0}$, $r_{j}$ is substituted for $r\geq r_{0}$, $\epsilon_{j}$ is substituted for $0$ and $\epsilon_{j}'$ is substituted for $\epsilon'(r)$. 
The constants in these inequalities depend on $A$ and $\rho$ only. 

Finally, observe that the continuous function $k(r)$, $r\geq r_{0}$, is positive nondecreasing, tending to $\infty$ as $r\to\infty$ (because $f$ is not a polynomial). In particular, we can determine its right inverse by the formula
\begin{equation}\label{e7471}
r(l):=\inf\{s\ : k(s)=l\},\ \ \ l\geq k(r_{0}):=k_{0}.
\end{equation}
Thus $r:[k_{0},\infty)\to [r_{0},\infty)$ is a continuous increasing function tending to $\infty$ as $k\to\infty$ and such that $k\circ r=id$. Substituting in the obtained inequalities and inequality (\ref{e746'}) $k$ instead of $k(r)$ and $r(k)$ instead of $r$ we obtain the required statements of Theorem \ref{te14} for $\rho<\infty$.
\\

(2) Assume now that $\phi_{f}(t):=m_{f}(e^{t})$ satisfies 
\begin{equation}\label{e748}
\lim_{t\to\infty}t^{2}\left(\frac{1}{\ln\phi_{f}(t)}\right)'=0.
\end{equation}

For each sufficiently large $v\in\Re$ by $s(v)\in\Re_{+}$ we denote a number such
that
$$
\frac{\phi_{f}(v)}{\phi_{f}(v-2s(v))}=e.
$$
Since $\phi_{f}$ is a continuous increasing function,
$$
s(v)=\frac{1}{2}\left(v-\phi_{f}^{-1}(\phi_{f}(v)/e)\right),
$$
i.e., $s(v)$ is a continuous in $v$ function.

Then similarly to (\ref{e731}) we obtain
\begin{equation}\label{e749}
\frac{1}{s(v)}\leq\frac{2e\phi_{f}'(v)}{\phi_{f}(v)}.
\end{equation}
From here arguing as in section 8.2 and using (\ref{e748}) for $\tilde s(v):=\min\{s(v),1\}$ we get 
\begin{equation}\label{e750}
\frac{1}{\tilde s(v)}\leq (m_{f}(e^{v}))^{\frac{\epsilon(v)}{v}},\ \ \ v\geq v_{0},
\end{equation}
where $\epsilon(v)$ is a positive continuous function in $v$ tending to $0$ as $v\to\infty$ and $v_{0}\in\Re$ is sufficiently large. 

Next, we determine continuous in $v$ functions 
$$
t(v):=e^{\tilde s(v)},\ \ \ \tilde r(v):=e^{v-\tilde s(v)}.
$$
Then as in section 8.2 we obtain
\begin{equation}\label{e751}
\ln\tilde r(v)\leq\epsilon'(v)[\ln m_{f}(t(v)\tilde r(v))]^{2},\ \ \ v\geq v_{0},
\end{equation}
where $\epsilon'(v)$ is a positive continuous function in $v$ such that $\lim_{v\to\infty}\epsilon'(v)=0$.
Also, in this case instead of Lemma \ref{below} by similar arguments we deduce that
\begin{equation}\label{e751'}
\ln\left(\frac{1+t(v)}{2\sqrt{t(v)}}\right)N_{f}(\tilde r(v),t(v))\geq (m_{f}(t(v)\tilde r(v)))^{1-\delta(v)}=:k(v),\ \ \ v\geq v_{0},
\end{equation}
for some positive continuous function $\delta(v)$ tending to $0$ as $v\to\infty$. Diminishing, if necessary, $1-\delta(v)$ we can assume that this function is increasing. In particular, the function $k(v)$, $v\geq v_{0}$, is increasing because $t(v)\tilde r(v)=e^{v}$.

From inequalities (\ref{e750}), (\ref{e751}), (\ref{e751'}) arguing as in section 8.2 
we obtain inequalities for $g\in{\cal F}_{p,q}(\tilde r(v);t(v);M_{f}(t(v)\tilde r(v)))$ with $p\leq k(v)$, $v\geq v_{0}$, similar to (\ref{e742}) in which $n_{j}^{1+\tilde\epsilon_{j}}$ is substituted for $(k(v))^{1+\tilde\epsilon(v)}$ 
for some continuous nonnegative function $\tilde\epsilon(v)$,
tending to $0$ as $v\to\infty$. As before, these inequalities give rise to inequalities of Theorem \ref{te13} with $\rho=\infty$ in which $n_{j}^{1+\epsilon_{j}}$ is substituted for $(k(v))^{1+\epsilon(v)}$, $v\geq v_{0}$, for some continuous nonnegative function $\epsilon(v)$ tending to $0$ as $v\to\infty$, and $r_{j}$ is substituted for $\frac{t(v)\tilde r(v)}{e^{2}}:=e^{v-2}$, $v\geq v_{0}$.

Finally, we write in the latter inequalities $k$ instead of $k(v)$. Then
since the function $k$ has a continuous increasing inverse
$s$, instead of $v$ and $e^{v-2}$ we obtain continuous increasing functions $s(k)$ and $r(k):=e^{s(k)-2}$, $k\geq k_{0}:=k(v_{0})$. This gives the required statements of the theorem (with $r_{0}:=r(k_{0})$). \ \ \ \ \ $\Box$\\
\\
{\bf 9.3. Proof of Corollary \ref{c14}.}
Inequality $\underline\tau(f)\leq 2$ follows directly from Theorem \ref{te13} (c) applied to restrictions to the graph of $f$ of polynomials $g$ of degrees $\lfloor n_{j}\rfloor$. (In this case $g\in {\cal F}_{p,q}(er_{j};e;M_{f}(e^{2}r_{j}))$ with $p=q=\lfloor n_{j}\rfloor\leq n_{j}$, $j\in\N$.)

Now, let us prove the lower bound for $\underline\tau(f)$.

By the definition the dimension of the space ${\cal P}_{k,n+1}$ of holomorphic polynomials of degree $k$ on $\Co^{n+1}$ is $d_{k,n+1}:=\frac{(n+k+1)!}{(n+1)!k!}$. Hence,
\begin{equation}\label{e752}
\lim_{k\to\infty}\frac{d_{k,n+1}}{k^{n+1}}=\frac{1}{(n+1)!}
\end{equation}
Let
$$
f(z):=\sum_{\alpha=0}^{\infty}f_{\alpha}z^{\alpha},\ \ \ 
z\in\Co^{n},
$$
be the Taylor series of $f$ at $0$. Here  $\alpha=(\alpha_{1},\dots,\alpha_{n})\in\Z_{+}^{n}$ and
$z^{\alpha}:=z_{1}^{\alpha_{1}}\cdots z_{n}^{\alpha_{n}}$.
The number of coefficients of the series at monomials of degrees $\leq p_{k}:=\left\lfloor\frac{k^{1+1/n}}{(n+2)^{1/n}}\right\rfloor$ is
$d_{p_{k},n}$. In particular,
\begin{equation}\label{e753}
\lim_{k\to\infty}\frac{d_{p_{k},n}}{k^{n+1}}=\frac{1}{(n+2)n!}
\end{equation}

Comparing (\ref{e752}) and (\ref{e753}) we conclude that
there is $k_{0}\in\N$ such that for each $k\geq k_{0}$ we have
$d_{p_{k},n}<d_{k,n+1}$.
Thus for each $k\geq k_{0}$ there is
$g\in {\cal P}_{k,n+1}$ such that the Taylor series of
$g_{f}(z):=g(z,f(z))$ has the form
\begin{equation}\label{e754}
g_{f}(z)=\sum_{\alpha=p_{k}+1}^{\infty}[g_{f}]_{\alpha}z^{\alpha},
\ \ \ \ z\in\Co^{n}.
\end{equation}
(Indeed, the linear map $\pi:{\cal P}_{k,n+1}\to{\cal P}_{p_{k},n}\ \!\!$, $\pi(g):=\sum_{\alpha=0}^{p_{k}}[g_{f}]_{\alpha}z^{\alpha}$, has a nonzero kernel because $d_{p_{k},n}<d_{k,n+1}$.) Also, since $f$ is not a polynomial, $g_{f}\not\equiv 0$.

Let $l\in {\cal L}_{n}$ be a complex line passing through $0$ such that 
$$
\sup_{l\cap\Bo^{n}}|g_{f}|=M_{g_{f}}(1).
$$
Let us identify $l$ with $\Co$. Then for the holomorphic function $h:=g_{f}|_{l}$ on $\Co$ by (\ref{e754}) we deduce that
the function
$$
\tilde h(z):=\frac{h(z)}{z^{p_{k}+1}}
$$
is holomorphic and $M_{\tilde h}(1)=M_{h}(1)\neq 0$. This implies
$$
M_{g_{f}}(1)=M_{\tilde h}(1)\leq M_{\tilde h}(e)\leq\frac{M_{h}(e)}{e^{p_{k}+1}}\leq \frac{M_{g_{f}}(e)}{e^{p_{k}+1}}.
$$
From here we have, see (\ref{e136}),
$$
m_{k}(e,f)\geq p_{k}+1.
$$
Therefore
$$
\liminf_{k\to\infty}\frac{m_{k}(e,f)}{k^{1+1/n}}\geq
\liminf_{k\to\infty}\frac{p_{k}+1}{k^{1+1/n}}=\frac{1}{(n+2)^{1/n}}>0.
$$
This shows that $\underline\tau(f)\geq 1+1/n$.\ \ \ \ \ $\Box$

\end{document}